\def \dd#1{{\bf#1}}

\def\cl#1{{\cal#1}}

%Simbolos matematicos.

\def\Min{\mathop{\rm Min}}

\def\Inf{\mathop{\rm Inf}}

%Macros.

\def\ouv#1{\smash{\mathop{#1}\limits^{\lower 1pt\hbox
{$\scriptscriptstyle\circ$}}}}

\def\hfl#1#2{\smash{\mathop{\hbox to 12mm{\rightarrowfill}}
\limits^{\scriptstyle#1}_{\scriptstyle#2}}}

%Titulos, enunciados.

\long\def\eno#1#2{\par\smallskip{\bf{#1}}{\it\ {#2}}\par\medskip}

\def\stit#1{\vskip 3mm plus 1mm minus 2mm {\bf{#1}}
		\smallskip}

\def\ref#1#2#3#4{{\bf #1}{\ #2}{\it ,\ #3}{,\ #4}\medskip}

%Dibujos

\def \picture #1 by #2 (#3){\midinsert \centerline 
{\vbox to #2{\hrule width #1 heigth 0pt 
depth 0pt \null \vfill \special {picture #3}}}\endinsert }

\def\scaledpicture #1 by #2 (#3 scaled #4) {{
\dimen0 =#1 \dimen1 =$2
\divide \dimen0 by 1000 \multiply \dimen0 by #4
\divide \dimen1 by 1000 \multiply \dimen1 by #4
\picture \dimen0 by \dimen1 (#3 scaled $4)}}

\def\figure #1 #2 #3 {\midinsert \vglue 3mm 
{\vbox to #3 {\hrule width 6cm height 0cm depth 0cm \vfill
{\special {picture #1 scaled #2}}}}\vglue 2mm \endinsert}

\magnification=1200

\input psfig.sty

\overfullrule=0pt

{\centerline {\bf ON TUBE-LOG RIEMANN SURFACES AND PRIMITIVES OF RATIONAL FUNCTIONS.}}

\bigskip
\bigskip

{\centerline {{\bf K. Biswas, R. P\'erez-Marco}}}

\bigskip
\bigskip

{\bf 1) Introduction.}

\medskip

The notion of {\it tube-log Riemann surface} was first introduced by R.Perez-Marco in order to solve several open
problems in holomorphic dynamics (see for example [PM1], [PM2], [PM3]). Somewhat informally, tube-log Riemann
surfaces may be described as Riemann surfaces constructed using a prescribed set of building blocks, including
complex planes $\dd C$ and complex cylinders $\dd C / 2\pi i \lambda \dd Z$ (hence the word ''tube''), by pasting
them together isometrically along 'slits', which are either finite line segments or half-lines. Each building
block comes equipped with a distinguished set of charts; for a complex plane the distinguished chart is the
identity $z \in \dd C \mapsto z$, while for a complex cylinder the distinguished charts are the locally defined
maps of the form $[z] \in \dd C / 2\pi i \lambda \dd Z \mapsto z$, and the changes of charts are translations. We
assume that the isometric pasting maps are translations in the distinguished charts. A tube-log Riemann surface
thus inherits a distinguished set of charts from its building blocks, for which the changes of charts are
translations. A tube-log Riemann surface also inherits a flat metric from its building blocks.

\medskip

Since any two distinguished charts differ by translations, any function $F : D \to \cl S$ taking values in a
tube-log Riemann surface $\cl S$ (where $D \subset \dd C$ is a planar domain) has a well-defined derivative
computed in the charts, $F' : D \to \dd C$. This allows us to write integral formulas for the uniformizations of
tube-log Riemann surfaces. The simplest example would be that of a single complex cylinder, for example $\cl S =
\dd C / 2\pi i \dd Z$. This is biholomorphic to the punctured plane $\dd C - \{ 0 \}$, with uniformization given
by the primitive $$\eqalign{ \dd C - \{ 0 \} & \to \dd C / 2\pi i \dd Z \cr
      z         & \mapsto \log z = \int_{1}^{z} {du \over u} \cr
}$$

\medskip

{\hfill {\centerline {\psfig {figure=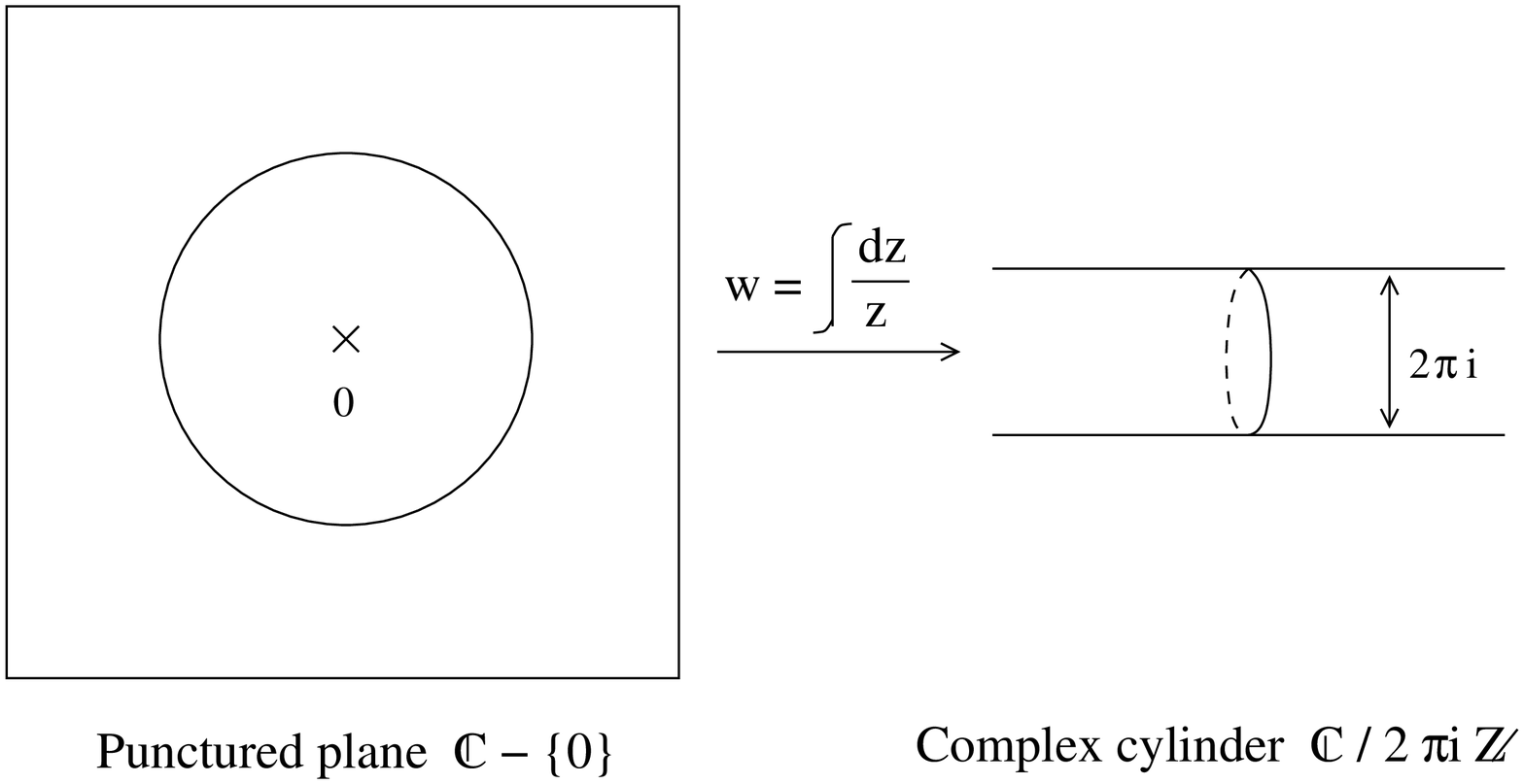,height=5cm}}}}

{\centerline {\bf Figure 1}}

\medskip

Another example is the tube-log Riemann surface $\cl S$ used by R.Perez-Marco in [PM1]. The surface $\cl S$ is
constructed from one complex cylinder and infinitely many complex planes, pasted together along half-line slits
(see figure below). It is also biholomorphic to the punctured plane $\dd C - \{ 0 \}$, with uniformization given
by the primitive $$\eqalign{ \dd C - \{ 0 \} & \to \cl S \cr
      z         & \mapsto \int_{1}^{z} {e^u \over u} \ du \cr
}$$

\medskip

{\hfill {\centerline {\psfig {figure=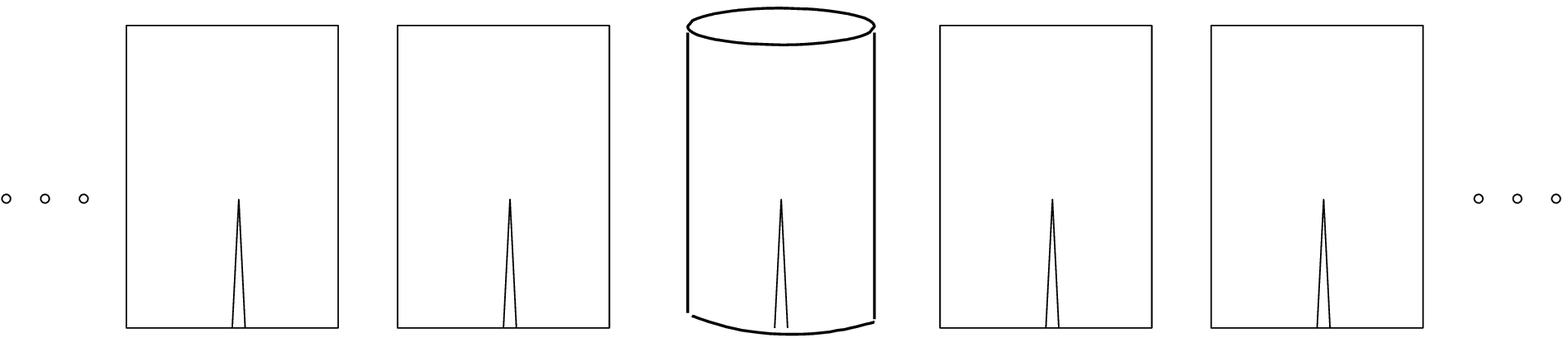,height=4cm,width=10cm}}}}

{\centerline {\bf Figure 2}}

\medskip

This geometry, and not another, proves the optimality of the diophantine condition ($(p_n/q_n)$ is the sequence of
convergents of the rotation number appearing in the problem) $$ \sum_{n=1}^{+\infty } {\log \log q_{n+1} \over
q_n} < +\infty \ , $$ in the Siegel problem of linearization of holomorphic dynamics with no strict periodic
orbits (see [PM1]).

\medskip

In this article we first define formally the notion of tube-log Riemann surface. We then study primitives of the
form $\int R(z) \ dz$, where $R(z)$ is a rational function of $z$. Our aim is to construct for each rational
function $R$ a tube-log Riemann surface $\cl S_R$ such that the uniformization of $\cl S_R$ is given by the
primitive $\int R(z) \ dz$. We show how to do this for a generic class of rational functions.

\medskip

A related study in the same spirit as this article which the reader may like to refer to, is that of {\it
Log-Riemann Surfaces} ([Bi-PM]). Log-Riemann surfaces form a subclass of tube-log Riemann surfaces, and are given
by pasting of only complex planes (they have no ''tubes'').

\bigskip

{\bf 2) Definition of tube-log Riemann surfaces.}

\medskip

{\bf 2.1) The basic building blocks.}

\medskip

We define here the building blocks allowed in the construction of tube-log Riemann surfaces.

\medskip

\eno{Definition (Building block).}{A building block is an object $B$ of one the five types (a)-(e) defined below.}

\medskip

{\bf (a) Complex cylinders:}

\medskip

Such a building block is of the form $B = \dd C / 2\pi i \lambda \dd Z$ for some $\lambda \in \dd Z$.

\medskip

{\bf (b) Complex half cylinders:}

\medskip

Such a building block is an open subset of a complex cylinder, of the form either $B = \{ [w] \in \dd C / 2\pi i
\lambda \dd Z : \hbox{ Im}(w/2\pi i \lambda) > 0 \}$ or $B = \{ [w] \in \dd C / 2\pi i \lambda \dd Z : \hbox{
Im}(w/2\pi i \lambda) < 0 \}$

\medskip

{\bf (c) Complex planes:}

\medskip

Such a building block is simply a copy of the complex plane, $B = \dd C$.

\medskip

The next two building blocks may seem a bit unusual, but they will arise naturally in the study of the primitives
of rational functions in the following sections. They require the notion of {\it log-Riemann surface}. For precise
definitions we refer the reader to [Bi-PM]; we give here a brief overview:

\medskip

A log-Riemann surface may be defined informally as a Riemann surface $\cl S$ with a distinguished set of charts
such that the trace of each chart is a slit complex plane, and such that any change of charts is the identity. The
charts on $\cl S$ thus paste together to give a globally defined local diffeomorphism $\pi : \cl S \to \dd C$,
which we refer to as the {\it projection mapping}. The log-Riemann surface $\cl S$ inherits a flat metric from its
charts, for which $\pi$ is a local isometry. We denote the completion of $\cl S$ with respect to this metric by
$\cl S^* = \cl S \sqcup \cl R$. We assume that the set $\cl R$ of points added is discrete, and call points of
$\cl R$ ramification points. It is shown in [Bi-PM] that for each ramification point $w^* \in \cl R$ there is a
small punctured ball $B(w^*, \epsilon) - \{w^*\}$ and an $n$, $2 \leq n \leq \infty$, such that $\pi$ maps $B(w^*,
\epsilon) - \{w^*\}$ to $B(\pi(w^*), \epsilon) - \{\pi(w^*)\}$ as an $n$-fold covering; we call $w^*$ a
ramification point of order $n$. Any two Euclidean segments meeting at $w^*$ delimit in this ball two angular
sectors (three if $n = \infty$) centered at $w^*$ of amplitudes $\theta, \phi >0$ such that $\theta + \phi =
2n\pi$. We can now define for log-Riemann surfaces the analogue of planar polygons:

\medskip

{\bf (d) Log-polygons:}

\medskip

Such a building block is an simply connected open subset of a log-Riemann surface $B \subset \cl S$ such that:

(i) Its boundary in $\cl S^*$, $\partial B \subset \cl S^*$, is homeomorphic to a circle.

(ii) For some $N \geq 3$, there are $N$ points on the boundary $\partial B$ such that the $N$ open arcs of
$\partial B$ delimited by them (each arc joins a point to the next in the circular ordering) are isometric to
finite open Euclidean segments.

The $N$ points are called the {\it vertices} and the $N$ arcs the {\it sides} of the log-polygon $B$. At each
vertex $w$ of $B$, the two sides of $B$ meeting at $w$ delimit, in a small ball around $w$, angular sectors
centered at $w$, of which one must be contained in $B$ (for a ball small enough). The amplitude of this sector is
called the {\it interior angle} at $w$.

\medskip

{\bf (e) Quotients of log-polygons:}

\medskip

Given a log polygon $P \subset \cl S$ with $N$ sides, one is allowed to construct such a building block $B$ as
follows:

If possible, take a strict subset of the sides of $P$, grouped in pairs $l_j, l'_j, j=1,\dots,k$ with $2k < N$,
such that for $j=1,\dots,k$, the Euclidean segments $\pi(l_j) \subset \dd C$ and $\pi(l'_j) \subset \dd C$ differ
by translations ($\pi : \cl S \to \dd C$ being the projection mapping of $\cl S$). Define $B$ to be the quotient
$$B = (P \cup (\cup_{j=1}^k l_j \cup l'_j)) /\sim ,$$ where we paste each $l_j$ to $l'_j$ by the map whose
expression in the distinguished charts on $\cl S$ is a translation.

We observe that if $B^*$ is the completion of $B$ with respect to the flat metric on $B$, then $B^* - B$ consists
of $N - 2k$ connected components, isometric to the $N - 2k$ sides of $\overline{P}$ not pasted together. We call
these components {\it boundary components} of the quotiented log-polygon $B$.

\bigskip

We remark that any building block is a Riemann surface equipped with a set of distinguished charts for which all
changes of charts are translations, and these charts induce a flat metric on the building block.

\medskip

{\bf 2.2) Cutting and pasting building blocks.}

\medskip

\eno{Definition (Cut or Slit).}{A cut or slit is a curve $\gamma$ in a building block $B$, $\gamma \subset B$,
which is isometric to either a closed half-line or a closed finite Euclidean segment.}

\medskip

\eno{Definition (Slit Building block).}{A slit building block is an open subset $B'$ of a building block $B$, of
the form $B' = B - \Gamma$, where $\Gamma = \sqcup_i \gamma_i$ is a disjoint union of cuts $\gamma_i$ in $B$ which
is locally finite, ie any compact subset of $B$ meets only finitely many cuts $\gamma_i$. The union of cuts
$\Gamma$ is allowed to be empty. We denote by $B'^*$ the completion of $B'$ with respect to the flat metric on
$B'$.}

Tube-log Riemann surfaces, which we imagine as being formed by cutting and pasting together building blocks, can
now be defined as follows:

\medskip

\eno{Definition (Tube-Log Riemann Surface).}{A tube-log Riemann surface is a Riemann surface $\cl S$ such that:

(1) $\cl S$ has a distinguished set of charts $\{ \phi_i : U_i \to \dd C \}_{i \in I}$ such that any change of
charts $\phi_i \circ \phi_j^{-1}$ is a translation.

(2) There is a collection $\{ V_k \}_{k \in J}$ of disjoint open subsets of $\cl S$ such that $$ \cl S =
\overline{\bigsqcup_{k \in J} V_k }, $$ and for each $V_k$ there is a slit building block $B'_k$ and a
biholomorphic map $\psi_k : B'_k \to V_k$ whose derivative computed in the distinguished charts on $B'_k$ and $\cl
S$ is identically equal to unity.

(3) We assume that each map $\psi_k$ extends to the completion of $B'_k$ minus a discrete set of points $D_k
\subset B'^*_k - B'_k$, to an injective mapping $\psi_k : B'^*_k - D_k \to \cl S$.
 }

\medskip

{\bf Remark.} The collection $\{ V_k \}_{k \in J}$ of slit building blocks of a tube-log Riemann surface $\cl S$
is not required to be unique; indeed in general it will not be unique, since some building blocks (for example
log-polygons) can be partitioned into smaller building blocks, and one can also change the directions of the
slits. We require however that at least one such decomposition exists. The condition (3) guarantees that all
slits, and boundaries of building blocks such as half-cylinders or log-polygons, are pasted together, except for
discrete sets of points on the slits or boundaries. The injectivity of the extended $\psi_k$ ensures that the two
'sides' of a slit are not identified in $\cl S$, so slits are non-trivial, and that there is no identification
within $\cl S$ of parts of the boundary of the same building block, for building blocks such as half-cylinders or
log-polygons (thus we exclude for example complex tori given by identifying opposite sides of a parallelogram).

\medskip

\eno{Definition.}{We regard two tube-log Riemann surfaces $\cl S_1$ and $\cl S_2$ as being equal if there is a
biholomorphic map $\psi : \cl S_1 \to \cl S_2$ whose derivative computed in the distinguished charts on $\cl S_2$
and $\cl S_2$ is identically equal to unity.}

\medskip

\bigskip

{\bf 3) Examples of tube-log Riemann surfaces for primitives of rational functions.}

\medskip

1. Let $R(z)$ have only one simple pole and no finite zeroes, say $R(z) = c/(z-a)$; this is essentially the same
as the example of the logarithm, with $\int R(z) \ dz = c \log (z-a)$ giving a uniformisation from the punctured
plane $\dd C - \{a\}$ to the complex cylinder $\dd C / 2\pi i c \dd Z$.

\medskip

2. Let $R(z)$ have two simple poles and no finite zeroes, say $R(z) = 1/(z-z_1)(z-z_2)$; splitting into
partial fractions and integrating we have
$$\eqalign{
\int {dz \over (z-z_1)(z-z_2)} & = {1 \over (z_1 - z_2)} [ \log (z - z_1) - \log(z - z_2) ] \cr
                               & = {1 \over (z_1 - z_2)} \log \left( { z-z_1 \over z-z_2 } \right). \cr
}$$
So $F(z) = \int R(z) dz$ is given by a Moebius transformation $z \mapsto (z - z_1)/(z - z_2)$, which maps the
punctured sphere $\overline{\dd C} - \{ z_1,z_2 \}$ to the punctured plane $\dd C - \{0\}$, composed with a
logarithm, which again maps to a complex cylinder.

\medskip

So though the degree of $R$ is $2$ in this example and $1$ in the first, both give the same tube-log Riemann
surfaces, namely complex cylinders. This is best explained by considering the 1-forms $R(z) \ dz$ instead of the
functions $R$. The rational function $R(z) = 1/z$ has only one pole at $z = 0$, so the 1-form $dz / z$ has a pole
there as well; however, near $z = \infty$, in terms of the coordinate $\xi = 1/z$, we have $dz / z = \xi \cdot
(-d\xi / \xi^2) = -d\xi / \xi$, so the form $dz/z$ also has an additional pole at $z = \infty$. Thus the 1-form
$dz/z$ has two poles and no zeroes; this is indeed also the case with the 1-form $dz/(z-z_1)(z-z_2)$, which has
two poles at $z=z_1$ and $z=z_2$, and no zeroes (the form is regular and non-zero near $z = \infty$), so the forms
$R(z) \ dz$ are of the same type in examples $1$ and $2$. In fact the Moebius transformation occuring above, $z
\mapsto (z - z_1)/(z - z_2)$, conjugates the forms $dz/(z-z_1)(z-z_2)$ and $dz/z$, taking the poles of the one to
those of the other.

\medskip

These remarks indicate that it is more appropriate to classify the primitives $\int R(z) \ dz$ according to the
types of the 1-forms $\omega = R(z) \ dz$ rather than those of the rational functions $R$. The 1-forms may be
broadly classified by the number $n$ of poles they have (counted with multiplicity); the number of zeroes is then
$(n-2)$ (counted with multiplicity). The geometry of the tube-log Riemann surfaces obtained will depend not only
on the total number of poles and zeroes, but also on their multiplicities.

\medskip

Before proceeding with more examples, it is worth carrying out a discussion which will be useful later of the
behaviour of an arbitrary 1-form near a simple pole:

\medskip

Let $\phi = f(z) \ dz = (\lambda/z + a_0 + a_1 z + a_2 z^2 + \dots) \ dz$ be a holomorphic 1-form with a
simple pole at $z = 0$. Then one can always make near $z=0$ an analytic change of variables $z = \xi +
O(\xi^2)$ to conjugate $\phi$ to the form $\lambda \ d\xi/\xi$, ie so that
$$
\phi = f(z) \ dz = \lambda \ d\xi/\xi.
$$
One way of seeing this is to 'rearrange' the above equation to derive the ODE that the change of variables
needs to satisfy,
$$
{dz \over d\xi} = {1 \over \xi} \cdot {\lambda \over f(z)} ,
$$
which it is easily seen has a unique analytic solution near $\xi=0$ with initial condition $z(0) = 0$.

\medskip

\eno{Definition (Pole petal).}{We define the {\it pole-petal} of the 1-form $\phi$ associated to the pole $z=0$ to
be the domain in the $z$-plane $P = z(\{ |\xi| < r \})$, where $r$ is the radius of convergence around $\xi = 0$
of the change of variables $z = z(\xi)$. Let $F(z) = \int f(z) \ dz$ be a primitive of $\phi$ near $z = 0$ (taking
values in $\dd C / 2\pi i \lambda \dd Z$); since the composition $\xi \mapsto z=z(\xi) \mapsto F(z(\xi)) = \lambda
\log \xi$ is a univalent function of $\xi$ in $\{ 0 < |\xi| < r \}$, so is $z = z(\xi)$. Thus the 'petal' $P$ is a
simply connected domain around the pole $z=0$ on which $\phi$ is conjugate to $\lambda \ d\xi/\xi$; moreover any
primitive $F$ of $\phi$ maps $P - \{0\}$ univalently to a half-cylinder $F(P) = \{ w \in \dd C / 2\pi i \lambda
\dd Z : \hbox{ Re} (w/\lambda) < \log r \} \subset \dd C / 2\pi i \lambda \dd Z$.}

\medskip

We return to the case of interest to us, ie that of a meromorphic 1-form $\omega = R(z) \ dz$ on the sphere; from
the above remarks, we see that the tube-log Riemann surface $\cl S$ of a primitive $\int R(z) \ dz$ should contain
half-cylinders, one for each simple pole; these should be somehow pasted together to form the surface $\cl S$. We
note that if for a pole the change of variables $z(\xi)$ is entire, so $R = \infty$, and its petal corresponds to
a full cylinder, then the form $\omega$ must have exactly two poles and be either as in example 1 or example 2;
from now on we exclude this trivial case, assuming that $\omega$ has a number of poles $n \geq 3$, so that $R <
\infty$ for every petal.

\medskip

\eno{Proposition.}{Let $\omega = R(z) \ dz$ be as above. We have:

(1) Pole-petals corresponding to distinct poles are disjoint.

(2) The boundary of a pole-petal contains no poles of $\omega$, but at least one zero of $\omega$.

(3) The boundary of a pole-petal is a piecewise analytic Jordan curve, analytic except at the zeroes of $\omega$.

(4) For two distinct pole-petals $P, P'$ with associated residues $\lambda, \lambda'$ respectively, if $\lambda$
and $\lambda'$ are not colinear then $\partial P$ and $\partial P'$ can meet only at zeroes of $\omega$.

 }

\medskip

\stit{Proof:} We first observe that a pole-petal contains no poles of $\omega$ other than the pole associated to
it.

\medskip

(1) : Let $P, P'$ be distinct pole-petals with associated poles $p,p'$ and residues $\lambda, \lambda'$
respectively. Suppose there is a $z_0 \in P \cap P'$. Let $\gamma$ be the closed equipotential curve of $P$
passing through $z_0$, say with periodic parametrization such that $R(\gamma(t)) \gamma'(t) = 2\pi i \lambda, t
\in \dd R$. If $\lambda, \lambda'$ are colinear, then $\gamma$ is an equipotential curve of $P'$ as well
(negatively oriented), and the two components of $\overline{\dd C} - \gamma$ containing $p,p'$ are contained in
$P$ and $P'$ respectively, so $\omega$ has only two poles $p,p'$, contradicting the hypothesis $n \geq 3$. If
$\lambda, \lambda'$ are not colinear, then we must have $\gamma(t) \to p'$ for either $t \to +\infty$ or $t \to
-\infty$, so, since $\gamma(t)$ is periodic, $p' \in \gamma \subset P$, contradicting the observation made above.

\medskip

(2) : That the boundary of a pole-petal contains no poles of $\omega$ follows from (1); this implies that the
boundary must contain at least one zero of $\omega$ (otherwise the map $\xi \mapsto z(\xi)$ could be extended
across every point of $\{ |\xi| = R \}$).

\medskip

(3) : Given (2), it is not hard to check that for each pole-petal $P$ of $\omega$, the change of variables $\xi
\mapsto z(\xi)$ extends to a continuous piecewise analytic map from the circle $\{ |\xi| = R \}$ to the boundary
$\partial P$ of the petal on the sphere, being analytic except for points which get mapped to zeroes of $\omega$;
so $\partial P$ is a piecewise analytic curve, possessing at each point $z \in \partial P$ (except for finitely
many 'corners') a tangent vector $X \in \dd C$ such that $R(z) \cdot X = 2\pi i \lambda$, where $\lambda$ is the
residue at the pole corresponding to $P$.

\medskip

(4) : The preceding remark implies that if two petals $P, P'$ have associated residues $\lambda, \lambda'$ which
are not colinear, then the boundaries $\partial P$ and $\partial P'$ can only intersect (if they do at all) at
their 'corners', ie at the zeroes of $\omega$ which lie on them (since at other points the boundaries are smooth
and non-tangential, so near an intersection point the domains themselves would have to intersect). $\diamondsuit$

\bigskip

We are now in a position to handle the following example:

\medskip

3. Let $\omega = (z-c_1) \ dz/(z - z_1)(z-z_2)(z-z_3)$, so that $\omega$ has $n = 3$ simple poles at
$z=z_1,z_2,z_3$, and a single zero at $z = c_1$. Let $P_1,P_2$ and $P_3$ be the pole-petals associated to the
poles $z_1,z_2$ and $z_3$. Assume that no two of the residues $\lambda_1,\lambda_2,\lambda_3$ at the three poles
are colinear; then by the remarks above, the boundaries of the petals meet precisely at the unique zero $z = c_1$.
Let $B$ be the complement in the sphere of the closed  connected set $\overline{P_1 \cup P_2 \cup P_3}$; $B$ is a
simply connected domain, so $\omega$ has a single-valued primitive $F(z) = \int R(z) \ dz$ in $B$. $F$ extends
analytically to all points of the boundary $\partial B = \partial P_1 \cup \partial P_2 \cup \partial P_3$,
mapping the boundary curves $\partial P_1, \partial P_2, \partial P_3$ to line segments corresponding to the
vectors $2\pi i \lambda_1, 2\pi i \lambda_2, 2\pi i \lambda_3$. Since $\lambda_1 + \lambda_2 + \lambda_3 = 0$, the
boundary $\partial B$ is in fact mapped to a triangle with these three vectors as sides; use of the argument
principle shows that $B$ is mapped univalently to the inside of this triangle.

\medskip

The primitive $F$ also extends analytically to each of the three petals, although it is not single-valued in them;
however, we can view it as mapping each petal $P_k$ to an open half-cylinder $C_k \subset \dd C / 2\pi i \lambda_k
\dd Z$. The tube-log Riemann surface $\cl S$ such that $F$ defines a uniformisation $F : \overline{\dd C} -
\{z_1,z_2,z_3,c_1\} \to \cl S$ can now be described as follows:

\medskip

For $k=1,2,3$ let $\overline{C_k}$ be the half-cylinder with boundary given by the closure in $\dd C / 2\pi i
\lambda_k \dd Z$ of $C_k$, so $\overline{C_k} = \{ w \in \dd C / 2\pi i \lambda_k \dd Z : \hbox{ Re}
(w/\lambda_k) \leq A_k \} \subset \dd C / 2\pi i \lambda_k \dd Z$ for some constant $A_k$.

\medskip

Let $T \subset \dd C$ be a closed triangle (ie including the interior and boundary) having for sides the
vectors $2\pi i \lambda_1, 2\pi i \lambda_2, 2\pi i \lambda_3$.

\medskip

Let $T' =  T - \{w_1,w_2,w_3\}$ be obtained from $T$ by deleting its vertices $w_1,w_2,w_3$. For $k=1,2,3$ let
$\overline{C_k}' = \overline{C_k} - \{p_k\}$ be obtained from $\overline{C_k}$ by deleting a single point $p_k$ on
its boundary. Then the boundaries of the $\overline{C_k}'$ s, which are isometric to open Euclidean line segments
given by the vectors $2\pi i \lambda_1, 2\pi i \lambda_2, 2\pi i \lambda_3$, correspond exactly to the boundary
segments of $T'$; we paste the boundary of each $\overline{C_k}'$ isometrically to the corresponding boundary
segment of $T'$, to obtain the tube-log Riemann surface $\cl S$.

\medskip

{\hfill {\centerline {\psfig {figure=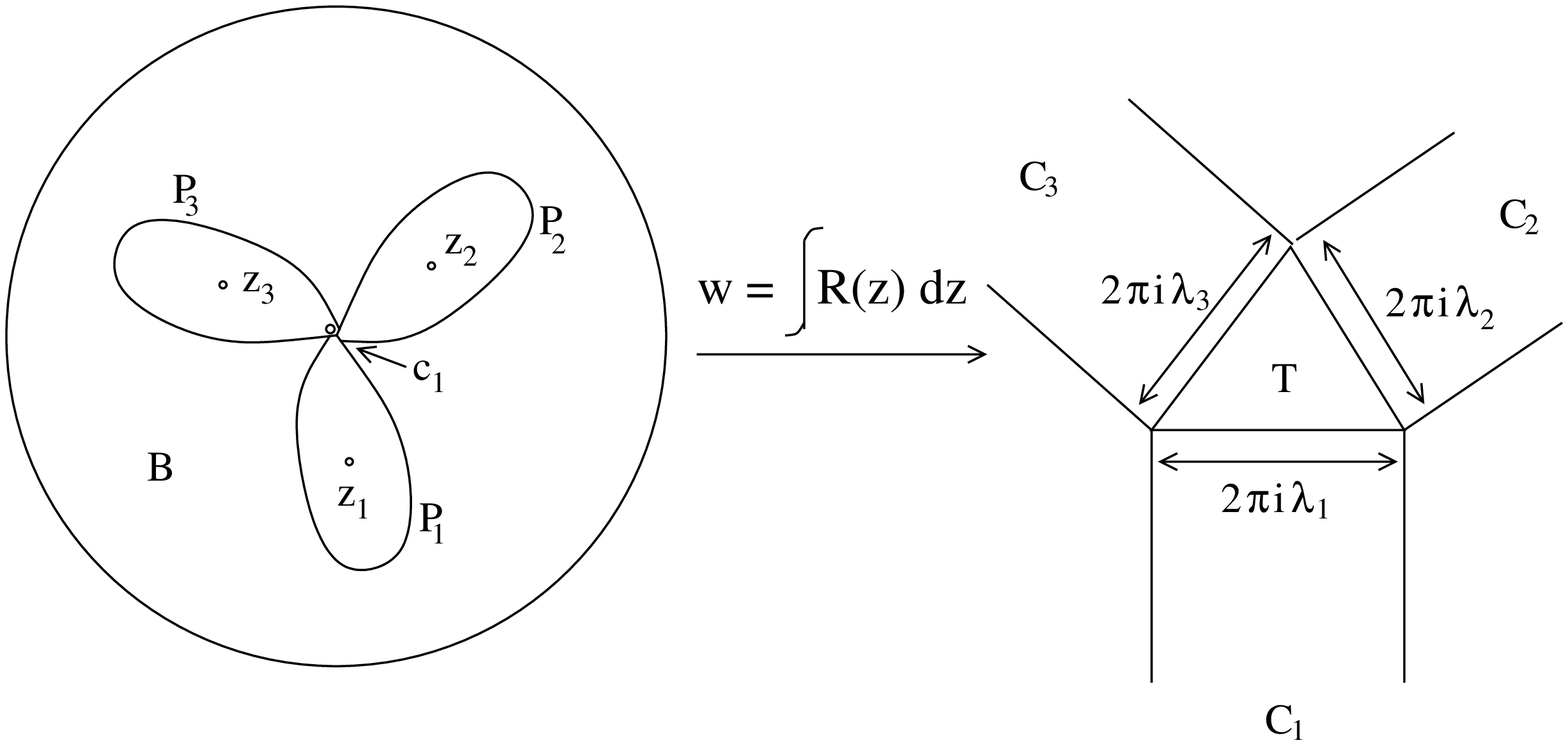,height=6cm}}}}

{\centerline {\bf Figure 3}}

\medskip

In this example we made the assumption that the residues $\lambda_k$ at the poles of $\omega$ were not colinear,
which is generically the case, and obtained a triangle with sides given by $2\pi i \lambda_k$, with each side
pasted to a half-cylinder. If the residues happen to be colinear, then the triangle degenerates into a 'triangle'
with angles $0,0$ and $\pi$, and the half-cylinders are pasted to each other along their boundaries. The following
example illustrates this case.

\medskip

4. Let $$ \omega = {2 \ dz \over z(z-1)(z+1)} = \left( {1 \over z-1} + {1 \over z+1} - {2 \over z} \right) \ dz $$
so $\omega$ has $3$ simple poles and one zero as in the preceding example. In this case however the residues, all
being real, are colinear, and instead of the boundaries of the petals $P_1,P_2,P_3$ meeting only at the critical
point, the petal associated to the pole $0$, say $P_1$, shares its boundary with the other two petals, ie
$\partial P_1 = \partial P_2 \cup \partial P_3$, and the petals give a partition of the sphere, $\overline{\dd C}
= \overline{P_1 \cup P_2 \cup P_3}$.

\medskip

To construct the associated tube-log Riemann surface $\cl S$ we first take as before the three half-cylinders with
boundary $\overline{C_1} \subset \dd C / 2\pi i \cdot (-2) \dd Z , \overline{C_2} \subset \dd C / 2\pi i \cdot (1)
\dd Z, \overline{C_3} \subset \dd C / 2\pi i \cdot (1) \dd Z$ corresponding to the petals $0,1$ and $-1$
respectively. As before we let $\overline{C_2}' = \overline{C_2} - \{p_2\}, \overline{C_3}' = \overline{C_3} -
\{p_3\}$ be obtained by deleting one point each from the boundaries of $\overline{C_2}, \overline{C_3}$. Let
$\overline{C_1}' = \overline{C_1} - \{q_1,q_2\}$ be given by deleting two points a distance $2\pi$ apart (for the
flat metric on $\dd C / 2\pi i \cdot (-2) \dd Z$) from the boundary of $\overline{C_1}$, so that the boundary of
$\overline{C_1}'$ consists of two disjoint line segments each of length $2\pi$. We paste the boundaries of
$\overline{C_2}'$ and $\overline{C_3}'$ isometrically to these two line segments, one to each, to obtain the
tube-log Riemann surface $\cl S$.

\medskip

In all the examples so far the associated tube-log Riemann surface $\cl S$ has been built up from half-cylinders
and polygons; in all cases the 1-form $\omega$ has had only one zero (not counting multiplicity). As we will see
shortly, the consideration of 1-forms with two or more distinct zeroes leads naturally to tube-log Riemann
surfaces with quotiented log-polygons as building blocks.

\bigskip

{\bf 4) Definition and uniqueness of the tube-log Riemann surface associated to the primitive $\int R(z) \ dz$.}

\medskip

\eno{Definition.}{ A tube-log Riemann surface $\cl S$ is said to be associated to the primitive $\int R(z) \ dz$,
where $R$ is a rational function, if there is a biholomorphic map $F : \overline{\dd C} - (Z \cup P) \to \cl S$
(where $Z,P$ are the zero and pole sets of $R$) such that its derivative $F' : \overline{\dd C} - (Z \cup P) \to
\dd C$ computed in the distinguished charts on $\cl S$ satisfies $$ F'(z) = R(z) $$}.

\medskip

\eno{Proposition.}{There is a unique tube-log Riemann surface associated to a primitive $\int R(z) \ dz$, ie if
$\cl S_1$ and $\cl S_2$ are two tube-log Riemann surfaces associated to the same primitive $\int R(z) \ dz$, then
$\cl S_1 = \cl S_2$ as tube-log Riemann surfaces.}

\medskip

\stit{Proof:} Given $\cl S_1, \cl S_2$ and the corresponding biholomorphic maps $F_1 : \overline{\dd C} - (Z \cup
P) \to \cl S_1, F_2 : \overline{\dd C} - (Z \cup P) \to \cl S_2$, since $F'_1(z) = R(z) = F'_2(z)$, it follows
that the biholomorphic map $F_2 \circ F^{-1}_1 : \cl S_1 \to \cl S_2$ has derivative in the distinguished charts
identically equal to unity. $\diamondsuit$

\medskip

We denote the unique tube-log Riemann surface associated to a primitive $\int R(z) \ dz$ by $\cl S_R$.

\bigskip

{\bf 5) Construction of the tube-log Riemann surface associated to $\int R(z) \ dz$.}

\medskip

Our problem is to try and construct, as in the examples in section 3, the tube-log Riemann surfaces $\cl S_R$ by
isometric pasting of building blocks for general primitives $\int R(z) \ dz$. One way of viewing the problem is
the following:

\medskip

The 1-form $R(z) \ dz$ gives, on the punctured sphere $\overline{\dd C} - (Z \cup P)$ a flat conformal metric
$|R(z)| \ |dz|$. This gives us a flat metric space, $(\overline{\dd C} - (Z \cup P), |R(z)| \ |dz|)$. The question
is, how can we describe this flat space? In particular, can this space be realized concretely by isometric pasting
of the building blocks defined in section 2?

\medskip

The general case of an arbitrary rational function $R$ however presents many complications, so we restrict
ourselves to describing a construction of the surfaces $\cl S_R$ that holds for a generic class of rational
functions.

\medskip

{\bf Notation.} Let $n \geq 4$ be the {\it degree} of the 1-form $R(z) dz$,
ie the number of poles (counted with multiplicity) of the 1-form. The number of zeroes
(counted with multiplicity) is then $n-2$. We denote the zero and pole sets
by $Z$ and $P$ respectively.

We denote the poles by $z_1, \dots, z_n$, the corresponding residues by
$\lambda_1, \dots, \lambda_n$ and the corresponding pole-petals by
$P_1, \dots, P_n$. We denote the zeroes by $c_1, \dots, c_{n-2}$.

\medskip

We restrict ourselves to 1-forms $R(z) dz$ satisfying the following
generic conditions:

\medskip

{\bf a) } All zeroes and poles are simple. There are thus $n$ distinct
poles and $n-2$ distinct zeroes.

\medskip

{\bf (b) } The sum of residues over any proper subset $I \subset \{1, \dots, n \},
\, I \not = \{1, \dots, n \}$ is non-zero,
$$
\sum_{i \in I} \lambda_i \not = 0
$$

\medskip

{\bf (c) } Any two such sums
$$
\sum_{i \in I} \lambda_i , \sum_{i \in J} \lambda_i
$$
over distinct subsets $I \not = J$ of $\{ 1,\dots, n \}$ are linearly independent
over $\dd R$.

\medskip

Our main result is the following :

\medskip

\eno{Main Theorem.}{ For the generic class described above, the tube-log Riemann surface $\cl S_R$ associated to a
primitive $\int R(z) \ dz$ is given by isometric pasting of $n$ half-cylinders $C_j \subset \dd C / 2\pi i
\lambda_j \dd Z, j=1,\dots,n$ to a quotiented log-polygon $B = P/\sim$ with $n$ boundary components. The
log-polygon $P$ has $3n-6$ sides and can be embedded in the log-Riemann surface $\cl S$ of a polynomial of degree
at most $2 \times (3n-6) + 1 = 6n - 11$.}

\medskip

By the {\it log-Riemann surface of a polynomial}, we mean a log-Riemann surface $\cl S$ whose uniformization is
given by a polynomial, ie there is a polynomial $Q \in \dd C[z]$ and a biholomorphic map $F : \dd C - \{ Q' = 0 \}
\to \cl S$ such that $\pi \circ F = Q$, where $\pi : \cl S \to \dd C$ is the projection mapping of $\cl S$.

\medskip

{\bf Terminology and tools.} We will talk of {\it geodesics} in the space $(\overline{\dd C} - (Z \cup P), |R(z)|
\ |dz|)$. By these we mean geodesics for the flat metric $|R(z)| \ |dz|$. We will also talk of geodesics $\gamma$
{\it in the direction of the vector} $\lambda \in \dd C$, by which we mean that the geodesic has velocity $$
R(\gamma(t)) \cdot \gamma'(t) = \lambda, $$ and geodesics {\it at an angle} $\theta$, meaning $$ \arg(
R(\gamma(t)) \cdot \gamma'(t) ) = \theta. $$ Given $\lambda \in \dd C$, these geodesics can also be thought of as
integral curves of the vector field $X_{\lambda}$ defined by $$ X_{\lambda}(z) := {\lambda \over R(z)} \ , \ z \in
\overline{\dd C} - (Z \cup P) $$ Note that these are curves which get mapped by any primitive $\int R(z) \ dz$ to
straight lines with direction vector $\lambda$.

\medskip

Before proceeding with the proof of the Theorem for general $n$, we first prove
the following Propositions and then work out the special case of degree $n = 4$.

\medskip

\eno{Proposition.}{For a 1-form $R(z) \ dz$ in the generic class described above, at
any zero $c_i$ there can be at most $2$ pole petals whose boundaries meet at
$c_i$.}

\medskip

\stit{Proof.} Let $c_i$ be a zero and suppose $k \geq 0$ pole petals meet at
$c_i$. Consider a local primitive $F(z) = \int R(z) \ dz$ defined near $c_i$,
and say $F(c_i) = 0$. $c_i$ is then a critical point of order 1 for $F$, which
maps the angle $2\pi$ at $c_i$ to an angle $4\pi$. Each petal boundary at $c_i$
is mapped to an angle $\pi$ by $F$, and these angles are disjoint; it follows
that $k \leq 4$. Moreover, if $k$ were equal to $4$, then the petal boundaries
would coincide and the corresponding residues would be parallel, contradicting
condition {\bf (c)} above. Hence $k \leq 3$. It remains to show that $k = 3$
cannot occur.

\medskip

Suppose then that $k = 3$, denote the petals by $P_1,P_2, P_3$ and the
corresponding residues by $\lambda_1, \lambda_2, \lambda_3$. We have the following
picture:

\medskip

{\hfill {\centerline {\psfig {figure=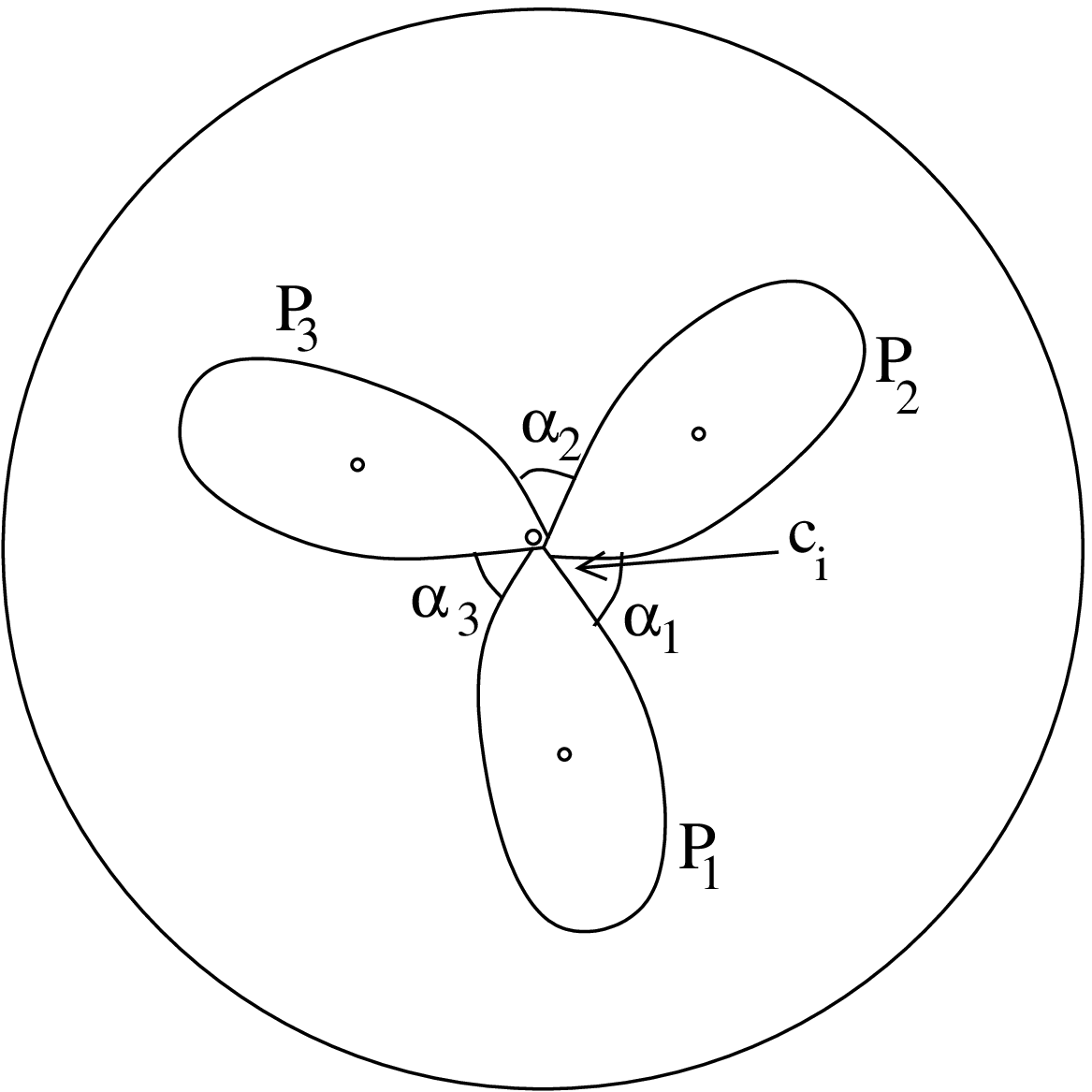,height=6cm}}}}

{\centerline {\bf Figure 4}}

\medskip

Since under the mapping $F$ the total angle at $c_i$ is equal to
$4\pi$ and the 3 petals contribute an angle $\pi$ each, it follows that
$$
2\alpha_1 + 2\alpha_2 + 2\alpha_3 = \pi
$$
(note angles double under $F$). In particular,
$2\alpha_1, 2\alpha_2, 2\alpha_3 < \pi$. Consider
small $\epsilon$-neighbourhoods (for the metric $|R(z)| |dz|$)
of the petal boundaries $\partial P_1, \partial P_2, \partial P_3$.
In the neighbourhood of $\partial P_j$, since $2\alpha_{j-1}, 2\alpha_{j} < \pi$ (the
indices are taken modulo 3 here), the integral curves to the vector field
$X_{2\pi i \lambda_j}$ starting from points in this neighbourhood exterior to $P_j$
must meet the boundaries $\partial P_{j-1}$ and $\partial P_{j+1}$. Moreover, if,
for 3 small constants $\epsilon_1, \epsilon_2, \epsilon_3 < \epsilon$,
$\gamma_j$ is the integral curve to $X_{2\pi i \lambda_j}$ starting from
a point at distance $\epsilon_j$ to $\partial P_j$, then the three
curves $\gamma_1, \gamma_2, \gamma_3$ meet to form a closed curve
$\gamma = \gamma_1' \cup \gamma_2' \cup \gamma_3'$, where $\gamma_j' \subset \gamma_j$
is a curve segment of $\gamma_j$.

\medskip

{\hfill {\centerline {\psfig {figure=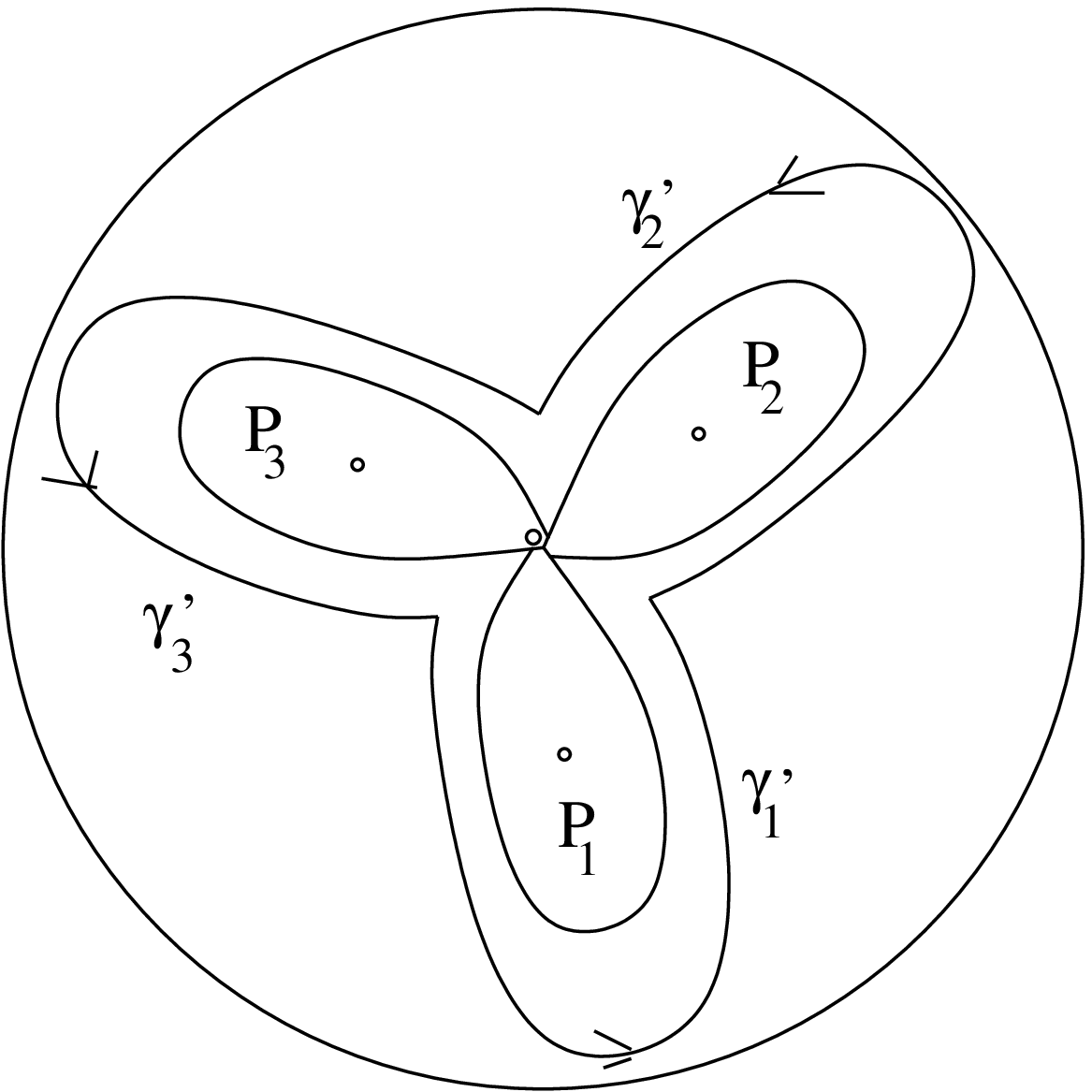,height=6cm}}}}

{\centerline {\bf Figure 5}}

\medskip

For each $\gamma_j'$ we have
$$
\int_{\gamma_1'} R(z) \ dz = (1 - \delta_j)2\pi i \lambda_j
$$
(where the $\delta_j$'s are small and depend on the $\epsilon_j$'s).
The $\epsilon_j$'s can be appropriately chosen so that
$\delta_1 = \delta_2 = \delta_3 = \delta > 0$ say.
In that case we have
$$\eqalign{
\int_{\gamma} R(z) \ dz & = \int_{\gamma_1'} R(z) \ dz + \int_{\gamma_2'} R(z) \ dz + \int_{\gamma_3'} R(z) \
dz \cr
 & = (1-\delta)2\pi i \lambda_1 + (1-\delta)2\pi i \lambda_2 + (1-\delta)2\pi i \lambda_3. \cr
}$$
On the other hand by Cauchy's residue formula we have
$$
\int_{\gamma} R(z) \ dz = 2\pi i \lambda_1 + 2\pi i \lambda_2 + 2\pi i \lambda_3.
$$
It follows from the above two equations that
$$
\lambda_1 + \lambda_2 + \lambda_3 = 0,
$$
a contradiction to condition {\bf (b)}. $\diamondsuit$

\medskip

Thus there can be at most $2$ pole-petals attached to a zero. In view of this
we make the following definitions:

\medskip

{\bf Definition.} We define $$\eqalign{ n_0 & = \hbox{ the number of zeroes of } R(z) \ dz \hbox{ with no petals
attached to them.} \cr n_1 & = \hbox{ the number of zeroes of } R(z) \ dz \hbox{ with one petal each attached to
them.} \cr n_2 & = \hbox{ the number of zeroes of } R(z) \ dz \hbox{ with two petals each attached to them.} \cr
}$$

\medskip

\eno{Proposition.}{For $n \geq 4$, all possible values of the triple $(n_0, n_1, n_2)$
are given by the triples
$$
(n_0 = j, \ n_1 = n-4-2j, \ n_2 = j+2) \ , \ 0 \leq j \leq [n/2] - 2
$$
}

\medskip

\stit{Proof.} Counting the number of zeroes and poles respectively of $R(z) \ dz$
gives the following two equations:
$$\eqalign{
n_0 + n_1 + n_2 & = n-2 \cr
n_1 + 2n_2 & = n \cr
}$$
Solving these two equations in three unknowns subject to the restrictions
$0 \leq n_0, n_1, n_2 \leq n-2$ gives the solutions listed above. $\diamondsuit$

\bigskip

{\bf 5.1) The case degree $n = 4$.}

\medskip

In this case, the above Proposition gives $$ n_0 = 0, \ n_1 = 0, n_2 = 2 $$ as the only possibility for the
arrangement of the pole-petals. Thus there are two pole-petals attached to each of the two zeroes, say $P_1, P_2$
to $c_1$ and $P_3, P_4$ to $c_2$.

\medskip

{\hfill {\centerline {\psfig {figure=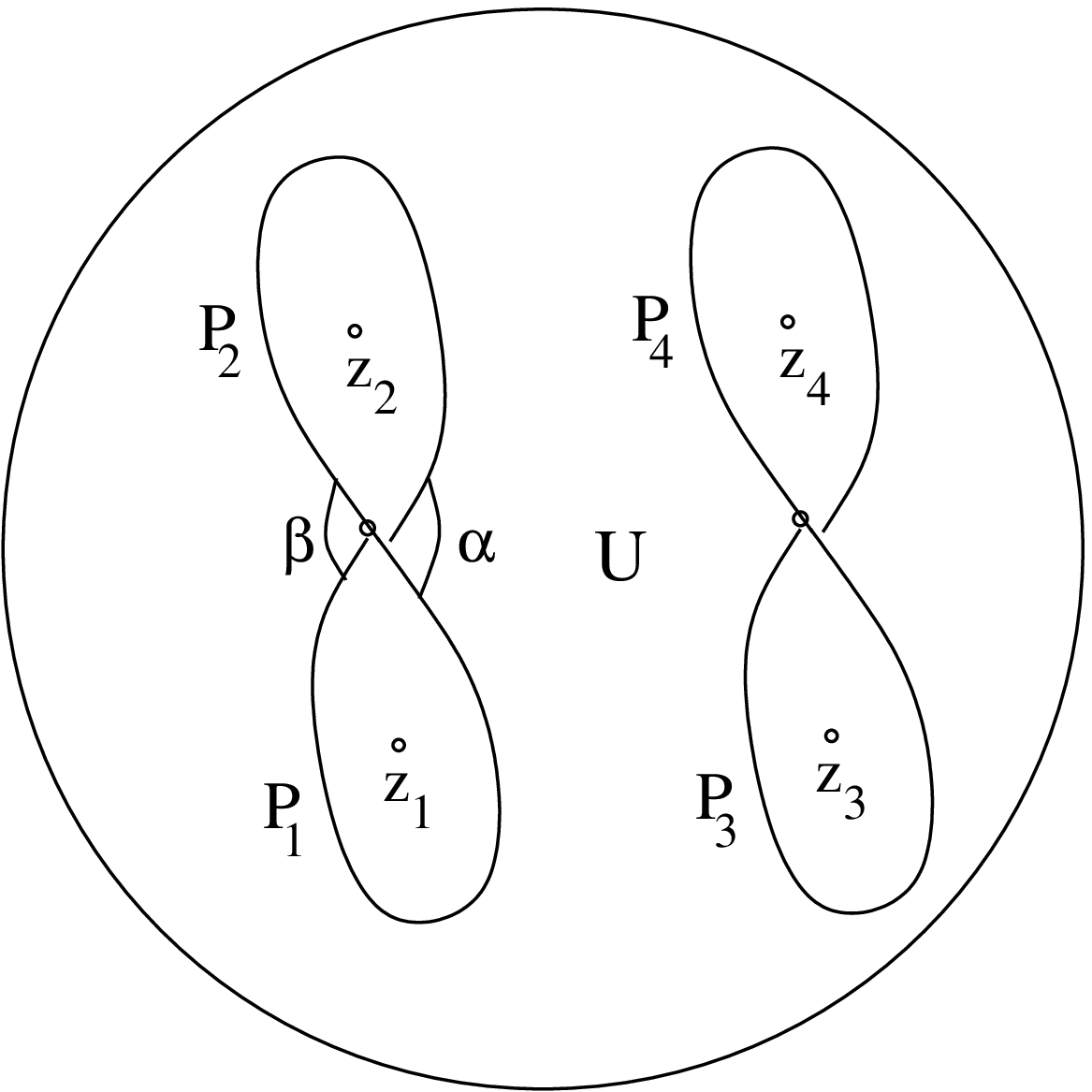,height=5cm}}}}

{\centerline {\bf Figure 6}}

\medskip

The annulus $U$ which separates these two pairs of pole-petals is of finite modulus; in fact it is isometric to
the quotient of a planar hexagon obtained by identifying a pair of parallel sides. More precisely, we have the
following Theorem:

\medskip

\eno{Theorem.}{Let $R(z) \ dz$ be a meromorphic 1-form of degree 4 on $\overline{\dd C}$ belonging to the generic
class defined in the previous section. There exists a planar hexagon $T$ with two pairs of sides, $2\pi i
\lambda_1, 2\pi i \lambda_2$ and $2\pi i \lambda_3, 2\pi i \lambda_4$, joined together by a pair of parallel sides
given by a vector $\mu \in \dd C$, such that the tube-log Riemann surface $\cl S_R$ is given by pasting
isometrically complex half-cylinders $C_j \subset \dd C/2\pi i \lambda_j \dd Z, j=1,\dots,4$ to the quotient of
$T$ given by identifying its pair of equal sides by the appropriate translation.}

\medskip

{\hfill {\centerline {\psfig {figure=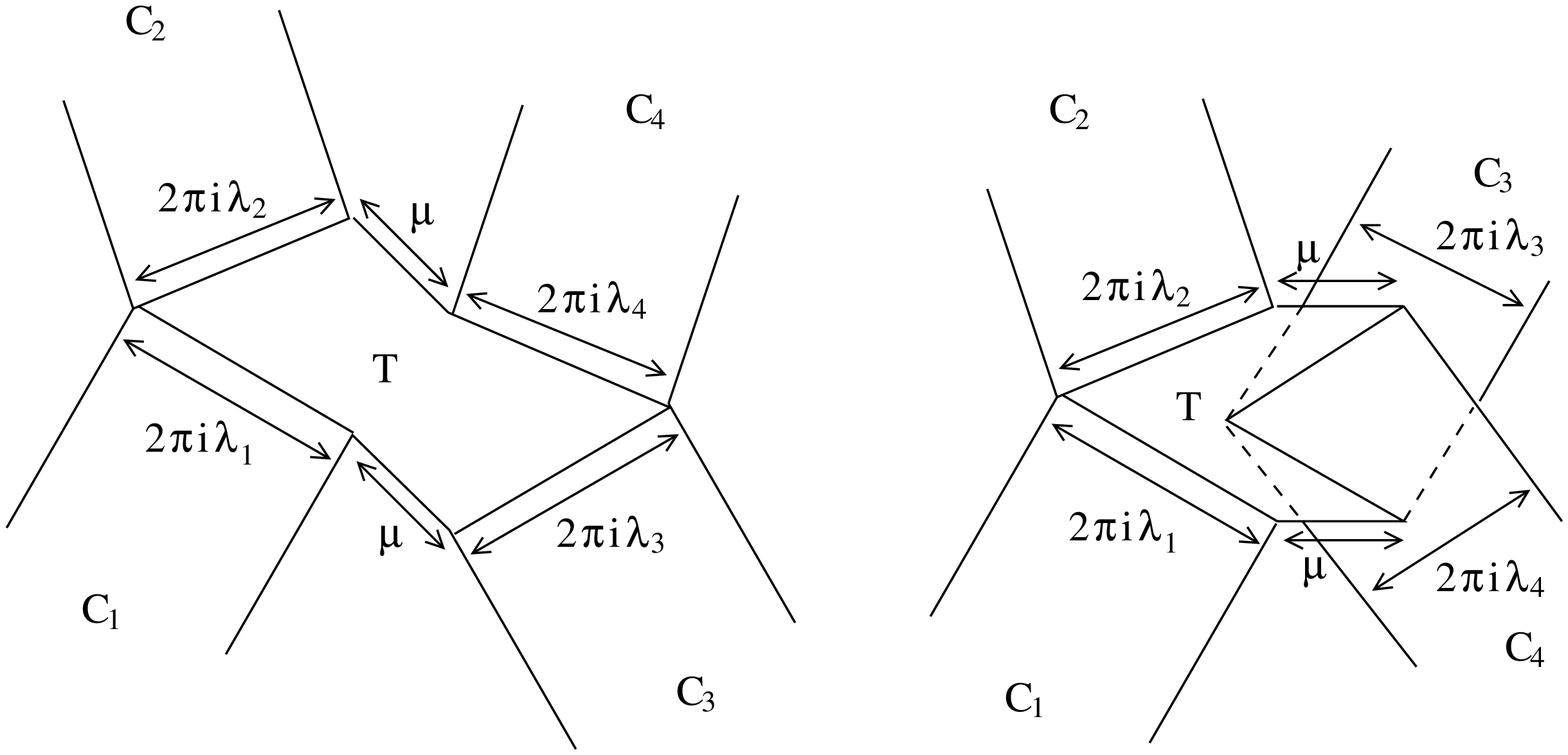,height=6cm}}}}

{\centerline {\bf Figure 7}}

\medskip

Thus in the special case $n = 4$, we can improve considerably on the general Theorem, in the sense that the
quotiented log-polygon given by the general Theorem can in fact be taken to be a quotiented planar polygon. The
above picture shows two distinct possible cases for the surface $\cl S_R$; in the figure to the left the whole
surface $\cl S_R$ can be realised isometrically as the quotient of a planar domain with boundary by translations,
whereas in the figure to the right this is not possible (the cylinders $C_3$ and $C_4$ 'overlap' when trying to
represent $\cl S_R$ isometrically as a planar quotient). Note that the hexagon $T$ need not be convex.

\medskip

\stit{Proof of Theorem.} The proof consists in 'drawing' parts of the surface $\cl S_R$ one by one.

\medskip

We already know that any primitive $\int R(z) \ dz$ maps each petal $P_j$
isometrically to a half-cylinder $C_j \subset 2\pi i \lambda_j$. It remains
to understand how it maps the annular region $U$ separating the two pairs
of petals.

\medskip

Considering the total angle at $c_1$ gives
$$
2\alpha + 2\beta = 2\pi,
$$
where $\alpha, \beta$ are the angles shown in figure 5. So one of the two
angles $2\alpha, 2\beta$ is less than $\pi$; suppose it is $2\alpha$,
so $2\alpha < \pi$. This implies that geodesics with direction
$2\pi i (\lambda_1 + \lambda_2)$ starting from points close to $c_1$
inside the angle $\alpha$ flow from one petal $P_1$ to the other $P_2$.
Thus we obtain a domain isometric to two half-cylinders $C_1, C_2$
joined to each other at an angle $2\alpha$ by pasting a small triangle
to their boundaries, as shown in the figure below.

\medskip

{\hfill {\centerline {\psfig {figure=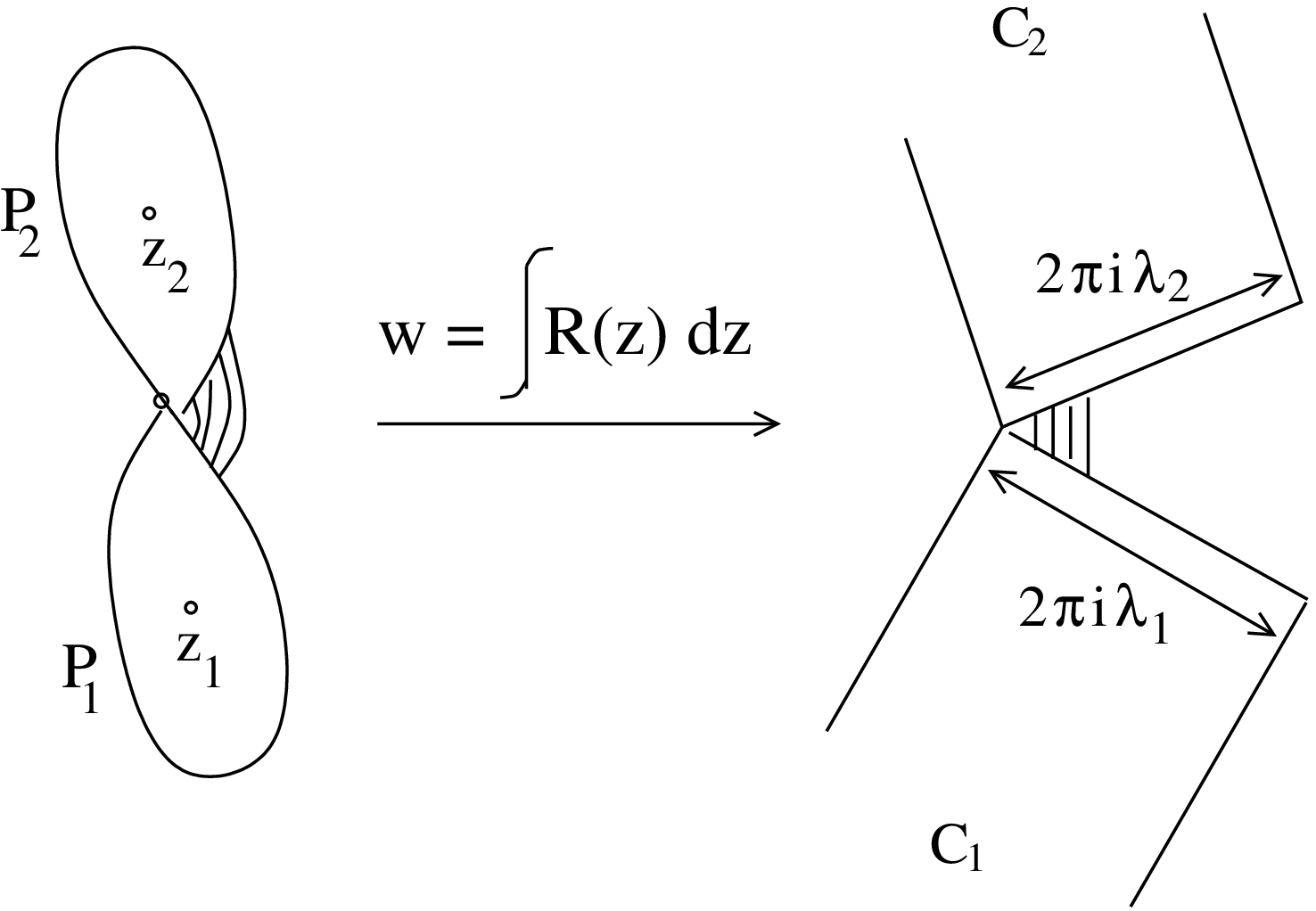,height=6cm}}}}

{\centerline {\bf Figure 8}}

\medskip

The only possible obstruction to continuing 'drawing' geodesics in this direction
from one petal boundary to the other and growing this triangle is posed by the zero
$c_2$, when one of these geodesics lands at $c_2$ before reaching the petal boundary
$\partial P_2$. We consider separately the two distinct cases that may occur:

\medskip

{\bf Case 1.} The geodesics can be continued without obstruction:

\medskip

In this case the geodesics can be continued till we arrive at a closed geodesic $\gamma_0$
entering and leaving the zero $c_1$ through the angle $\beta$. The region bounded by
$\gamma_0$ and the petal boundaries $\partial P_1, \partial P_2$ is isometric to a triangle
$T_1$ with sides $2\pi i \lambda_1, 2\pi i \lambda_2$ and $-2\pi i (\lambda_1 + \lambda_2)$,
and the whole region enclosed by $\gamma_0$ is isometric to two half-cylinders $C_1, C_2$
joined to each other at an angle $2\alpha$ by pasting this triangle to their boundaries.

\medskip

The geodesics with the same direction $2\pi i (\lambda_1 + \lambda_2)$ starting from
points close to $\gamma_0$ and to the exterior of the region bounded by $\gamma_0$ must
be smooth closed geodesics. This family of geodesics can be grown till we encounter
the second zero $c_2$; when this occurs we obtain a closed geodesic $\gamma_1$
starting and ending at $c_2$. The region bounded by $\gamma_0$ and $\gamma_1$ is made
up of smooth closed geodesics, and is isometric to a finite cylinder $C$ of the
form
$$
C = \{ w \in \dd C / 2\pi i (\lambda_1 + \lambda_2) \dd Z : a < \hbox{ Re} (w/(\lambda_1 + \lambda_2)) < A \}
\subset \dd C / 2\pi i (\lambda_1 + \lambda_2) \dd Z
$$
for some constants $a, A \in \dd R$.  Thus the half-cylinders $C_1, C_2$ are
pasted to the triangle $T_1$ whose free boundary side is pasted to
a finite cylinder $C$, as shown in the figure below.

\medskip

{\hfill {\centerline {\psfig {figure=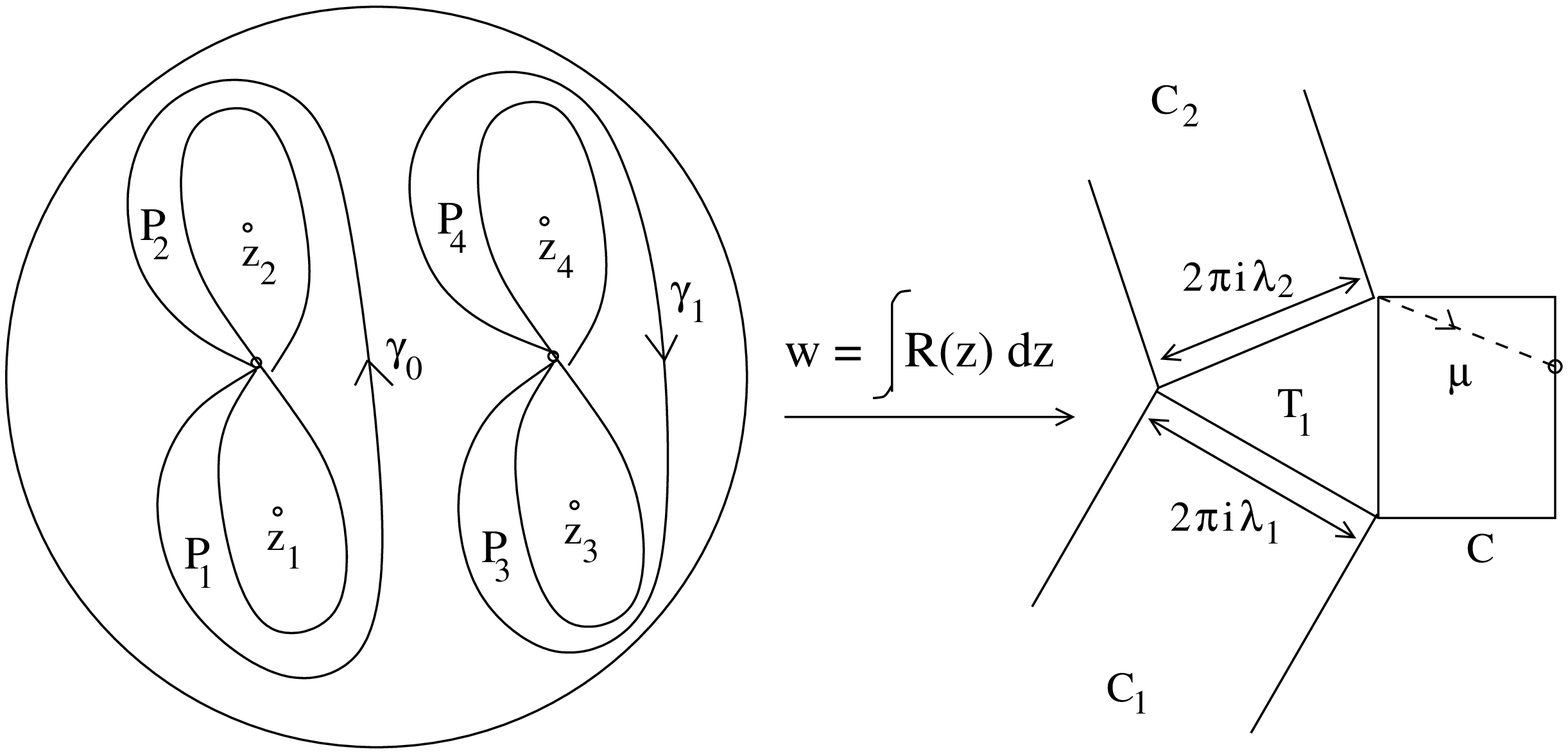,height=6cm}}}}

{\centerline {\bf Figure 9}}

\medskip

Moreover, the points $c_1, c_2$ correspond to points on the two boundaries
of the finite cylinder $C$ and can hence be joined by a geodesic isometric
to a line segment $\mu \in \dd C$ as shown in the figure. The finite cylinder
$C$, shown in the figure as the quotient of a rectangle, can also be
represented isometrically as the quotient of a parallelogram $L$ with one pair of
parallel sides equal to $\mu$ and the other to $2\pi i (\lambda_1 + \lambda_2)$.

\medskip

The only region that remains to be understood is the region bounded
by $\gamma_1$ and the petal boundaries $\partial P_3, \partial P_4$. This is a
simply connected domain whose boundary segments $\partial P_3, \partial P_4,
\gamma_1$ are mapped by any primitive $F(z) = \int R(z) \ dz$ to line segments
$2\pi i \lambda_3, 2\pi i \lambda_4$ and $2\pi i (\lambda_1 + \lambda_2) =
- 2\pi i (\lambda_3 + \lambda_4)$ respectively, and hence to the boundary of
a triangle $T_2$ with these three sides. It follows that this domain is mapped to
the interior of $T_2$; finally, the petals $P_3, P_4$ are mapped to
half-cylinders $C_3, C_4$ pasted isometrically along their boundaries to the
corresponding sides of $T_2$.

\medskip

Thus the tube-log Riemann surface $\cl S_R$ is given by pasting the two pairs of half-cylinders $C_1, C_2$ and
$C_3, C_4$ to the triangles $T_1$ and $T_2$ respectively, pasting $T_1$ and $T_2$ to the parallelogram $L$ and
finally identifying the sides of $L$ equal to $\mu$.

\medskip

{\hfill {\centerline {\psfig {figure=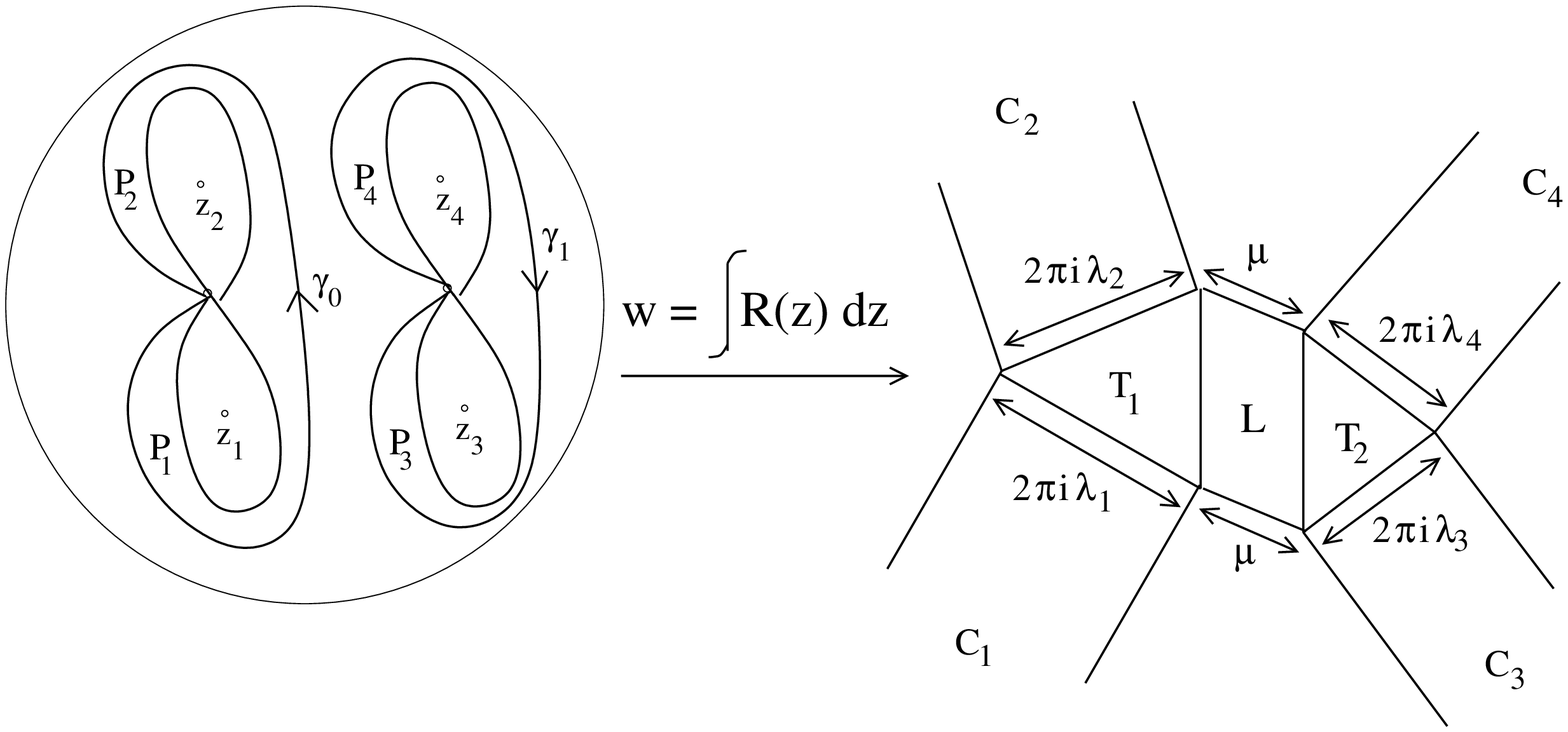,height=6cm}}}}

{\centerline {\bf Figure 10}}

\medskip

The hexagon $T$ of the Theorem is given by the union of the triangles
$T_1, T_2$ and the parallelogram $L$, $T = T_1 \cup L \cup T_2$, pasting
the appropriate sides of $T_1, T_2$ to those of $L$ isometrically.

\medskip

{\bf Case 2.} The geodesics with direction $2\pi i (\lambda_1 + \lambda_2)$
starting from points on $\partial P_1$ hit the zero $c_2$:

\medskip

In this case we obtain a domain containing the petals $P_1, P_2$ that
is isometric to two half-cylinders $C_1, C_2$ joined to each other by
pasting a triangle $T_1$ to part of their boundaries. The free boundary of the
triangle corresponds to two geodesics $\alpha_1, \alpha_2$
which start from points on $\partial P_1, \partial P_2$ and meet at
the zero $c_2$, as shown in the figure below:

\medskip

{\hfill {\centerline {\psfig {figure=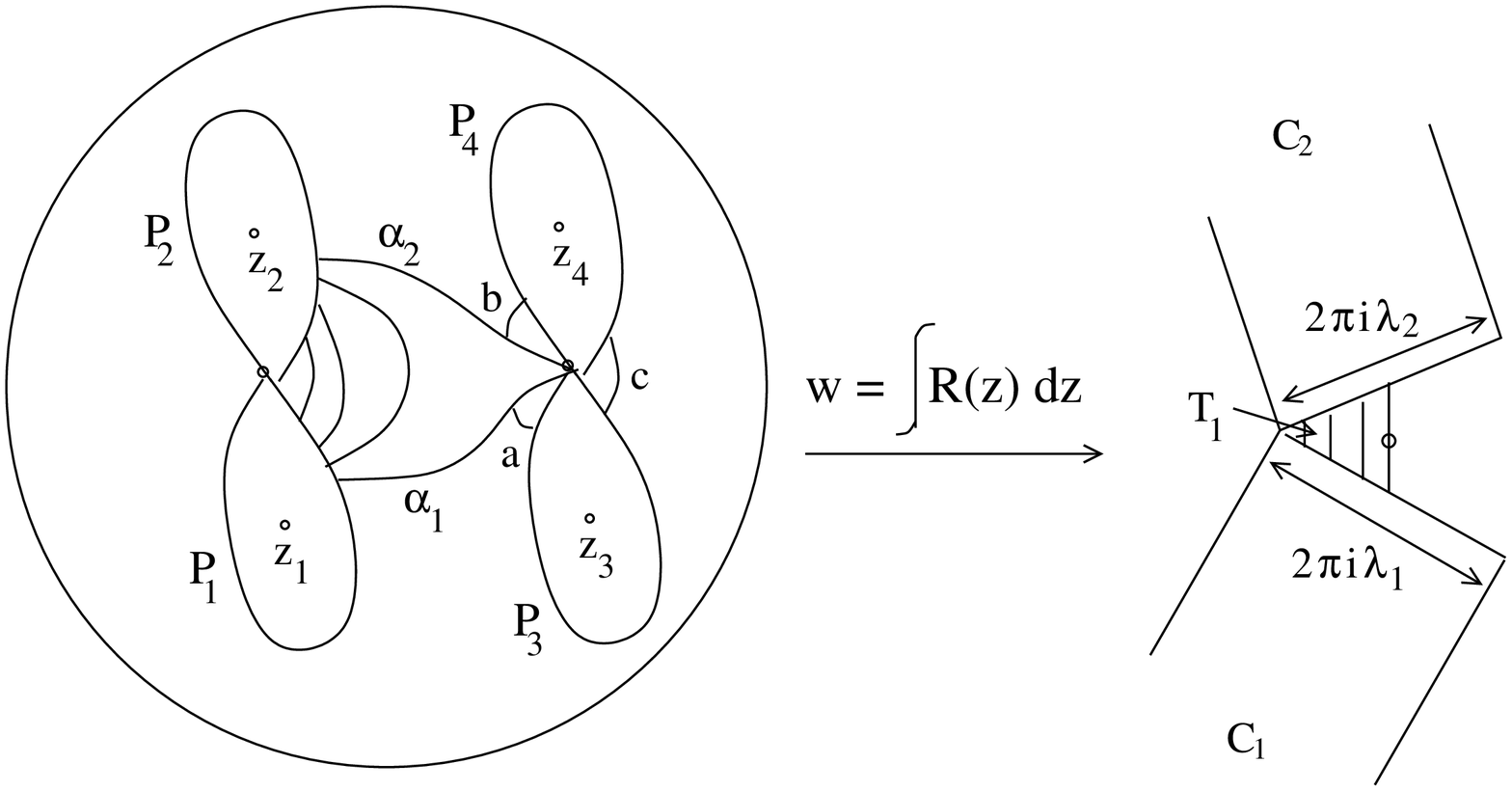,height=6cm}}}}

{\centerline {\bf Figure 11}}

\medskip

We observe that $2a + 2b + 2c = \pi$, so in particular $2a, 2b < \pi$.
Geodesics starting from the zero $c_2$ and leaving through the angle $a$
in a direction almost parallel to that of $\alpha_1$ must
meet $\partial P_1$, similarly geodesics leaving through the angle $b$
in a direction almost parallel to that of $\alpha_2$ must
meet $\partial P_2$. We can keep increasing the angles between these
geodesics and $\alpha_1, \alpha_2$ till we obtain geodesics $\beta_1,
\beta_2$ that start from $c_2$ and meet at $c_1$. This gives two domains
isometric to two triangles $T_2, T_3$ pasted to the triangle $T_1$ and the
half-cylinders $C_1, C_2$:

\medskip

{\hfill {\centerline {\psfig {figure=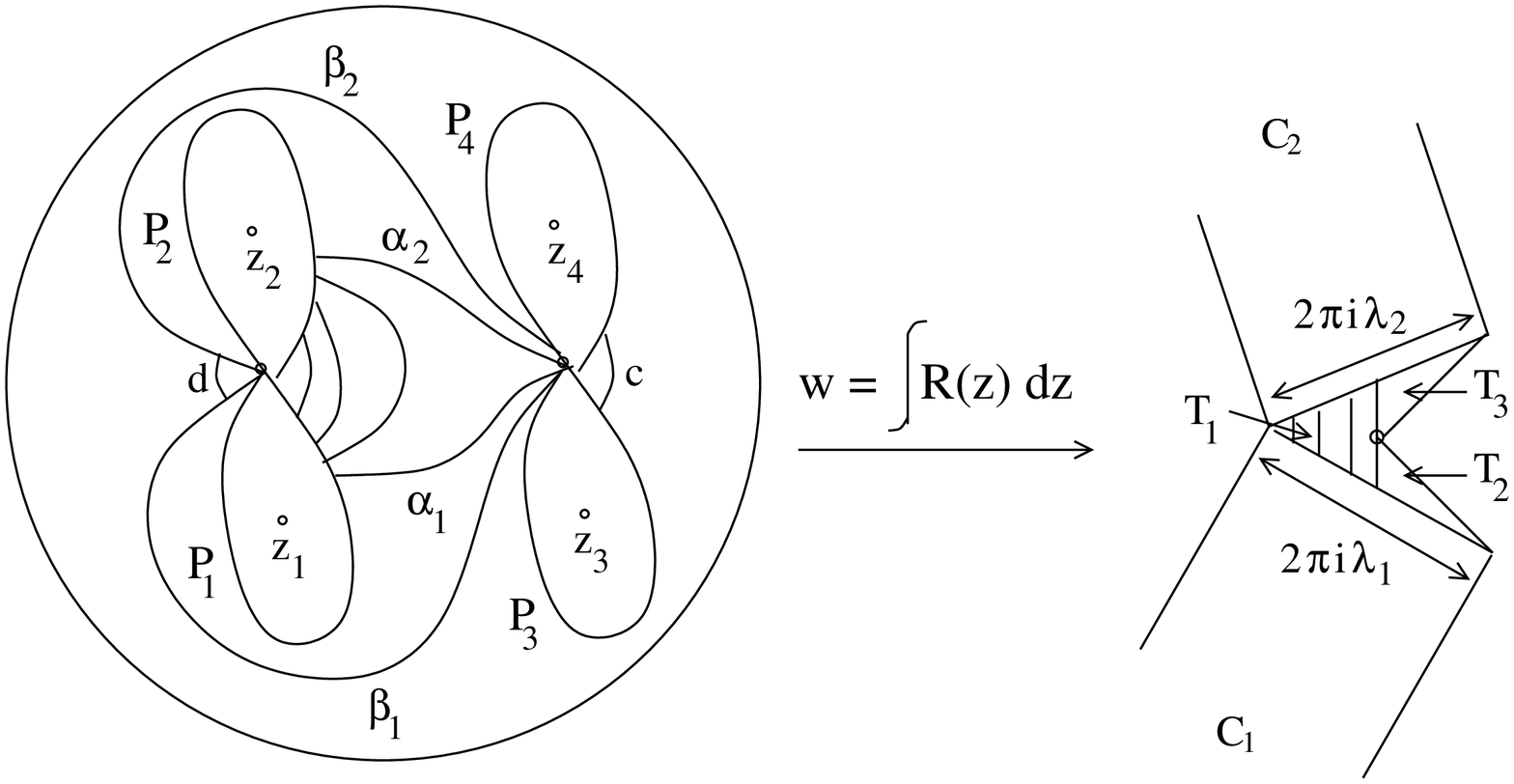,height=6cm}}}}

{\centerline {\bf Figure 12}}

\medskip

Since the angle $2a$ is less than $\pi$, now geodesics starting from
the zero $c_1$ and leaving through the angle $d$
in a direction almost parallel to that of $\beta_1$ must
meet $\partial P_3$. We can keep increasing the angles between these
geodesics and $\beta_1$ till we obtain a geodesic $\tau$ that starts from
$c_1$ and ends at $c_2$. This gives two simply connected domains $D_1, D_2$,
bounded by the curves $\tau, \beta_1, \partial P_3$ and by
$\tau, \beta_2, \partial P_4$ respectively. Since any primitive $\int R(z) \ dz$
defined in these domains maps the boundaries to boundaries of triangles, these
domains are isometric to two triangles $T_4, T_5$ pasted to the triangles $T_2,
T_3$ and the half-cylinders $C_3, C_4$ as shown in the figure below
(here $\mu \in \dd C$ is given by $\mu = \int_{\tau} R(z) \ dz$):

\medskip

{\hfill {\centerline {\psfig {figure=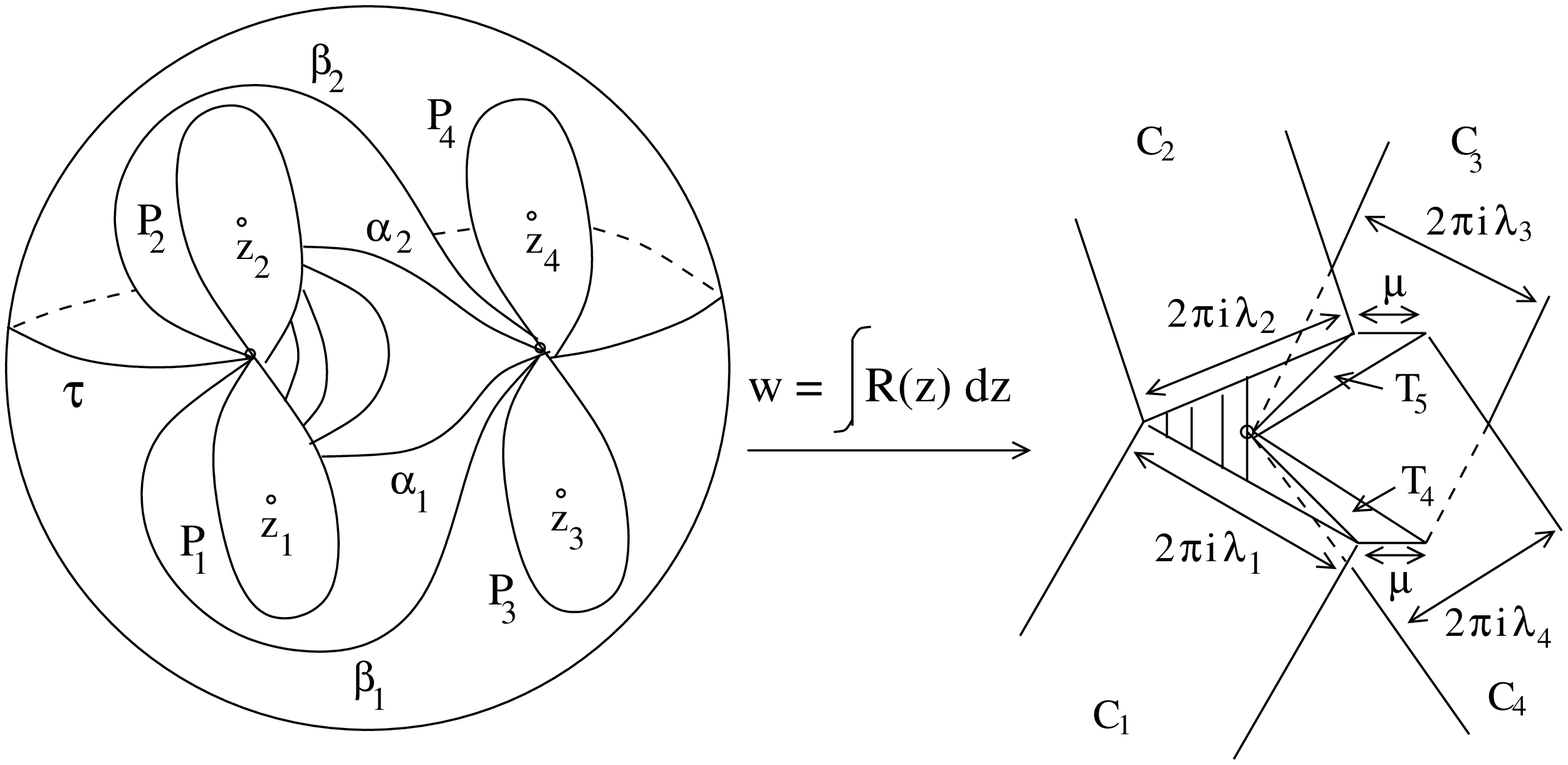,height=6cm}}}}

{\centerline {\bf Figure 13}}

\medskip

The hexagon $T$ of the Theorem is given by the union of the triangles
$T_1, \dots, T_5$, pasted isometrically along corresponding boundaries.
This ends the proof of the Theorem. $\diamondsuit$

\bigskip

{\bf 5.2) The general case of degree $n \geq 4$.}

\medskip

We now consider the general case of degree $n \geq 4$. It will be useful
to consider the region in the sphere $\overline{\dd C}$ exterior to the
open pole-petals $P_1, \dots, P_n$ as a metric space with the metric induced
by the infinitesimal metric $|R(z)| |dz|$.

\medskip

{\bf 5.2.1) The metric space $(X,d)$.}

\medskip

{\bf Definition.} We define the metric space $(X,d)$ by:
$$
X := \overline{\dd C} - (P_1 \cup \dots \cup P_n),
$$
and
$$
d(z_1, z_2) := \Inf_{\gamma} \ l(\gamma) \ , \ z_1, z_2 \in X
$$
where the infimum is taken over all rectifiable paths $\gamma$ in $X$
joining $z_1$ to $z_2$ and $l(\gamma)$ denotes the length of $\gamma$
computed with respect to the metric $|R(z)||dz|$.

\medskip

The space $X$ is compact, and the classical argument using
Ascoli-Arzela's Theorem applies in this setting to give

\medskip

\eno{Proposition.}{The distance between any two points of
$X$ is attained by a curve in $X$, ie given $z_1, z_2  \in X$
there exists a path $\gamma$ in $X$ joining $z_1$ to $z_2$
such that
$$
d(z_1,z_2) = l(\gamma).
$$
We call $\gamma$ a minimizing geodesic joining $z_1$ and $z_2$.}

\medskip

Since the metric $|R(z)| |dz|$ is flat away from the critical
points $Z$ where $R(z) dz$ vanishes, minimizing geodesics must
be isometric to Euclidean segments away from the critical points.
Minimizing geodesics can thus be described as follows:

\medskip

\eno{Proposition.}{Any minimizing geodesic is isometric either to a single
Euclidean segment or to a polygonal line which is a union of Euclidean
segments, with vertices at critical points.}

\bigskip

{\bf 5.2.2) Connecting the critical points by geodesics in $X$.}

\medskip

The space $X$ is not simply connected, since the connected components
of its complement in $\overline{\dd C}$ consist of distinct groups
of pole-petals attached to critical points, and it is not hard to see
from Proposition [] that there must be at least two such groups. We aim
to make 'cuts' in $X$ in such a way that the resulting region $D$ is
simply connected, so that we will be able to define a single-valued primitive
$\int R(z) \ dz$ in $D$.

\medskip

We will choose $(n-3)$ 'cuts' $\gamma_1, \dots, \gamma_{n-3}$ in order to
connect the $(n-2)$ critical points $c_1, \dots, c_{n-2}$. Each 'cut'
$\gamma_k$ will be a geodesic in $X$ joining a pair of critical points.
The geodesics $\gamma_k, k=1,\dots,(n-3)$, will be chosen inductively
according to the following algorithm:

\medskip

{\bf Step 1. \ At stage $k=1$:}

\medskip

Choose critical points $c^{(1)},c^{(2)} \in Z$ such that
$$
d(c^{(1)},c^{(2)}) = \Min_{1 \leq i \neq j \leq n-2} d(c_i, c_j)
$$
and choose $\gamma_1$ to be a minimizing geodesic in $X$ joining $c^{(1)}$ to
$c^{(2)}$ such that
$$
d(c^{(1)},c^{(2)}) = l(\gamma_1).
$$
If $k+1=2=n-2$ (ie if $n=4$), then we stop at this point, otherwise
we set $k = k+1$ and proceed as follows.

\medskip

{\bf Step 2. \ At stage $k \geq 2$:}

\medskip

We assume that $k$ critical points $c^{(1)},c^{(2)},\dots,c^{(k)} \in Z$
and $(k-1)$ geodesics $\gamma_1,\dots,\gamma_{k-1}$ have been chosen.

\medskip

We choose a critical point $c^{(k+1)} \in Z - \{ \ c^{(1)},c^{(2)},\dots,c^{(k)} \ \}$,
distinct from those already chosen, such that it is closest to those
already chosen, in the sense that
$$
d(c^{(k+1)}, \{ \ c^{(1)},c^{(2)},\dots,c^{(k)} \ \} ) = \Min_{c \in Z - \{ \ c^{(1)},c^{(2)},\dots,c^{(k)} \
\} } d(c, \{ \ c^{(1)},c^{(2)},\dots,c^{(k)} \ \} ).
$$
(by the distance $d(c^{(k+1)}, \{ \ c^{(1)},c^{(2)},\dots,c^{(k)} \ \} )$ we
mean the distance between a point $x$ and a compact set $A$, defined as usual by
$d(x,A) = \Min_{y \in A} d(x,y)$).

\medskip

Let $c^{(i_k)} \in \{ \ c^{(1)},c^{(2)},\dots,c^{(k)} \ \}$ be a critical point
such that
$$
d(c^{(k+1)}, \{ \ c^{(1)},c^{(2)},\dots,c^{(k)} \ \} ) = d(c^{(k+1)}, c^{(i_k)} )
$$
and choose $\gamma_k$ to be a minimizing geodesic joining $c^{(k+1)}$ to $c^{(i_k)}$ such
that
$$
d(c^{(k+1)}, c^{(i_k)} ) = l(\gamma_k).
$$
We observe that $\gamma_k$ cannot pass through any critical points other
than $c^{(k+1)}$ and $c^{(i_k)}$ (this would clearly contradict the
definition of $c^{(k+1)}, c^{(i_k)}$), and is thus isometric to a single
Euclidean segment.

\medskip

If $k+1 = n-2$, then all critical points $c^{(1)},\dots,c^{(n-2)}$ and all
the required geodesics $\gamma_1, \dots, \gamma_{n-3}$ have been chosen,
so we stop at this point, otherwise we set $k = k+1$ and repeat Step 2.

\bigskip

The above algorithm, when it terminates, gives us $(n-3)$ geodesics
$\gamma_1, \dots, \gamma_{n-3}$. We can show that these geodesics do
not intersect (except possibly at their endpoints):

\medskip

\eno{Proposition.}{ Assume $n \geq 5$ (so that $n-3 \geq 2$ and
at least two geodesics $\gamma_1, \gamma_2$ have been chosen). Then
for any two distinct geodesics $\gamma_k, \gamma_{k'}, 1 \leq k \neq k' \leq n-3$
either
$$
\gamma_k \cap \gamma_{k'} = \phi
$$
or
$$
\gamma_k \cap \gamma_{k'} = \{ c \}
$$
where $c \in Z$ is a critical point which is a common endpoint of $\gamma_k$ and
$\gamma_{k'}$.
}

\medskip

\stit{Proof:} Suppose that two such geodesics $\gamma_k, \gamma_{k'}, k \neq k',$
meet at a point $z_0$ which is not a critical point,
$$
z_0 \in \gamma_k \cap \gamma_{k'} \ , \ z_0 \notin Z.
$$
We may assume $k' > k$. The endpoints of $\gamma_k$ are critical points
$c^{(i_k)}, c^{(k+1)}$ while those of $\gamma_{k'}$ are critical points
$c^{(i_{k'})}, c^{(k'+1)}$.
Removing the point $z_0$ from the curves $\gamma_k,
\gamma_{k'}$ disconnects each into two connected components, one containing
each endpoint. We let
$$\displaylines{
\gamma_{k,1} := \hbox{Connected component of } \gamma_k - \{z_0\} \hbox{ containing } c^{(i_k)} \cr
\gamma_{k,2} := \hbox{Connected component of } \gamma_k - \{z_0\} \hbox{ containing } c^{(k+1)} \cr
}$$
and similarly
$$\displaylines{
\gamma_{k',1} := \hbox{Connected component of } \gamma_{k'} - \{z_0\} \hbox{ containing } c^{(i_{k'})} \cr
\gamma_{k',2} := \hbox{Connected component of } \gamma_{k'} - \{z_0\} \hbox{ containing } c^{(k'+1)} \cr
}$$
The situation is as shown in the figure below.

\medskip

{\hfill {\centerline {\psfig {figure=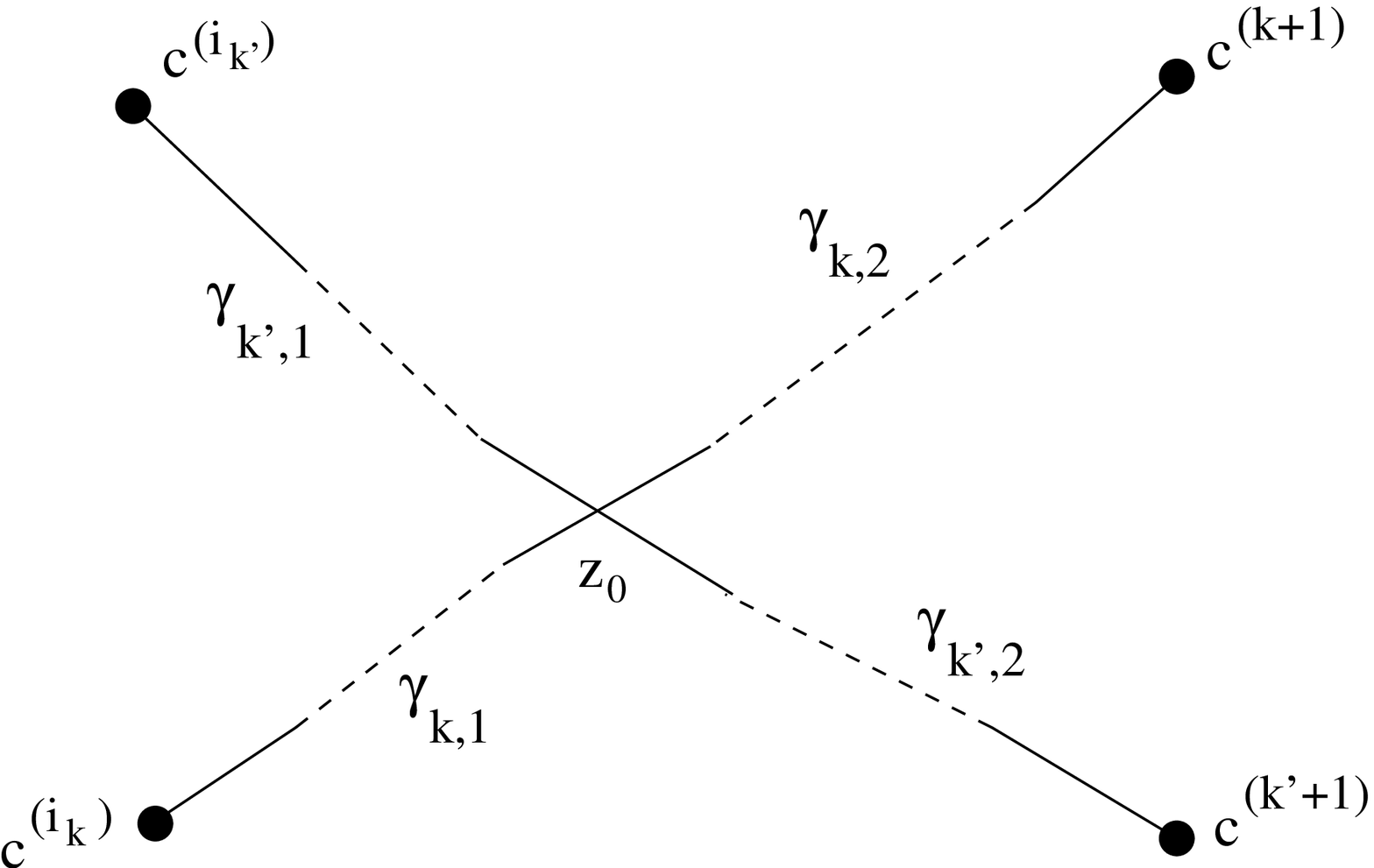,height=6cm}}}}

{\centerline {\bf Figure 14}}

\medskip

Observing that $c^{(i_k)} \in \{ \ c^{(1)}, \dots, c^{(k)} \ \} \subset \{ \ c^{(1)}, \dots, c^{(k')} \ \}$,
we have
$$\eqalign{
l(\gamma_{k',1})+l(\gamma_{k',2}) = d(c^{(i_{k'})}, c^{(k'+1)}) &= \Min_{c \in \{ \ c^{(1)}, \dots, c^{(k')}
\ \} } d(c, c^{(k'+1)}) \cr
 & \leq d(c^{(i_k)}, c^{(k'+1)}) \cr
 & \leq l(\gamma_{k,1})+l(\gamma_{k',2}) \cr
}$$
and hence
$$
l(\gamma_{k',1}) \leq l(\gamma_{k,1}).
$$
The curves $\gamma_{k',1}$ and $\gamma_{k,2}$ are geodesic segments meeting at $z_0$; since $z_0$ is not a
critical point, we can modify the curve $\gamma_{k',1} \cup \gamma_{k,2}$ slightly in a
neighbourhood of $z_0$ to construct a curve $\gamma$ joining $c^{(i_{k'})}$ to
$c^{(k+1)}$ such that
$$
l(\gamma) < l(\gamma_{k',1} \cup \gamma_{k,2}).
$$
Therefore
$$\eqalign{
d(c^{(i_{k'})}, c^{(k+1)}) \leq l(\gamma) & < l(\gamma_{k',1})+l(\gamma_{k,2}) \cr
          & \leq l(\gamma_{k,1}) + l(\gamma_{k,2}) \cr
          & = d( c^{(i_k)}, c^{(k+1)} ) \cr
          & = \Min_{ c \in \{ \ c^{(1)}, \dots, c^{(k)} \ \} } d(c, c^{(k+1)}) \cr
}$$
It follows that $c^{(i_{k'})} \notin \{ \ c^{(1)}, \dots, c^{(k)} \ \}$.
This implies that
$$\eqalign{
l(\gamma_{k,2})+l(\gamma_{k,1}) = d(c^{(k+1)}, c^{(i_k)}) & = \Min_{c \in Z - \{ \ c^{(1)}, \dots, c^{(k)} \
\} } d( c, \{ \ c^{(1)}, \dots, c^{(k)} \ \} ) \cr
                                                          & \leq d( c^{(i_{k'})}, \{ \ c^{(1)}, \dots,
c^{(k)} \ \} ) \cr
                                                          & \leq d( c^{(i_{k'})}, c^{(i_k)} ) \cr
                                                          & \leq l(\gamma_{k',1}) + l(\gamma_{k,1}) \cr
}$$
and hence
$$
l(\gamma_{k,2}) \leq l(\gamma_{k',1}).
$$
As before, we can modify the curve $\gamma_{k,2} \cup \gamma_{k',2}$ slightly in a
neighbourhood of $z_0$ to construct a curve $\gamma'$ joining $c^{(k+1)}$ to
$c^{(k'+1)}$ such that
$$
l(\gamma') < l(\gamma_{k,2} \cup \gamma_{k',2}).
$$
This implies that
$$\eqalign{
d( c^{(k+1)}, c^{(k'+1)} ) \leq l(\gamma') & < l(\gamma_{k,2}) + l(\gamma_{k',2}) \cr
                                           & \leq l(\gamma_{k',1}) + l(\gamma_{k',2}) \cr
                                           & = d( c^{(i_{k'})}, c^{(k'+1)} ) \cr
                                           & = \Min_{ c \in \{ \ c^{(1)}, \dots, c^{(k')} \ \} } d( c,
c^{(k'+1)} ) \, , \cr
}$$
which is a contradiction, since $c^{(k+1)} \in \{ \ c^{(1)}, \dots, c^{(k')} \ \}$. $\diamondsuit$

\bigskip

{\bf 5.2.3) Construction of the log-polygon $P$.}

\medskip

We now define the open set $$ D := X - \left( \bigcup_{j=1}^{n} \partial P_j \cup \bigcup_{k=1}^{n-3} \gamma_k
\right) $$. Its not hard to see from the preceding Proposition that $D$ is connected. Since the complement of $D$
is clearly connected (the $\gamma_k$'s connect all the critical points and hence all the petals $P_k$ as well), it
follows that $D$ is simply connected. The figure below shows an example of the domain $D$ for degree $n = 6$.

\medskip

{\hfill {\centerline {\psfig {figure=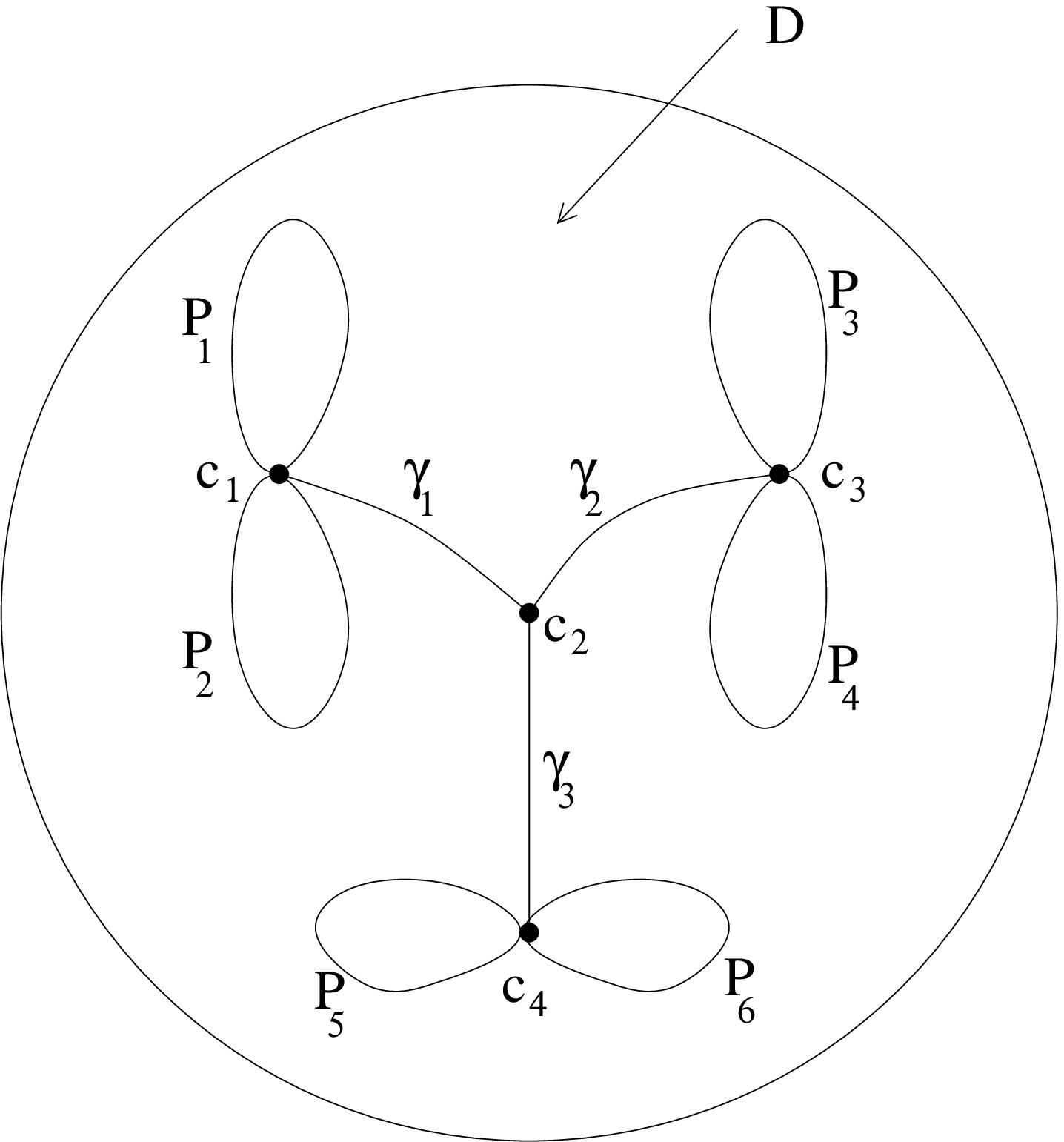,height=6cm}}}}

{\centerline {\bf Figure 15}}

\medskip

Fixing a base point $z_0 \in D$ we can define in $D$ a single-valued primitive $F : D \to \dd C$ given by $$ F(z)
:= \int_{z_0}^{z} R(t) \ dt \ , \ z \in D. $$

\medskip

The infinitesimal metric $|R(z)| |dz|$ on $D$ induces a global metric $d_1$ on $D$ defined by $$ d_1(z_1, z_2) :=
\Inf_{\gamma} \ l(\gamma) \ , \ z_1, z_2 \in D $$ where the infimum is taken over all rectifiable paths $\gamma$
in $D$ joining $z_1$ to $z_2$ and $l(\gamma)$ denotes the length of $\gamma$ computed with respect to the metric
$|R(z)||dz|$. Let $D^* = D \sqcup \cl E$ denote the Caratheodory compactification of $D$, given by adding the
circle of prime ends $\cl E$ of $D$. The metric $d_1 : D \times D \to \dd R$ extends in a unique way to a metric
$d_1 : D^* \times D^* \to \dd R$. The map $F : D \to \dd C$, which is a local isometry of $(D, d_1)$, extends to a
local isometry of $(D^*, d_1)$, which we continue to denote by $F : D^* \to \dd C$.

\medskip

Since the boundary $\partial D$ of $D$ is locally connected, prime ends of $D$ correspond to points of $\partial
D$ (not necessarily in a one-to-one fashion). For each point of a petal boundary $\partial P_j$ there is one prime
end that corresponds to it (except for the critical points), while for each point of a curve $\gamma_k$ there are
two prime ends that correspond to it (again excepting the critical points). The circle of prime ends $\cl E$ can
be written as a finite union of arcs intersecting only at endpoints, each of which is mapped isometrically by $F$
to a compact Euclidean segment. To each petal boundary $\partial P_j$ corresponds one such arc, mapped by $F$ to a
Euclidean segment with endpoints differing by the vector $2\pi\lambda_j \in \dd C$, while to each curve $\gamma_k$
correspond two such arcs, mapped by $F$ to a Euclidean segment with endpoints differing by a vector $2\pi\tau_k
\in \dd C$. There are thus $n + 2(n-3) = 3n - 6$ such arcs.

\medskip

\eno{Theorem.}{There exists a log-Riemann surface $\cl S$ and a log-polygon $P \subset \cl S$, such that $D^*$
embeds isometrically into $\cl S$, via an isometry $\tilde{F} : D^* \to P \cup \partial P \subset \cl S^*$ mapping
$D$ onto $P$ and $\cl E$ onto $\partial P$. The log-Riemann surface $\cl S$ can be taken to be the log-Riemann
surface of a polynomial of degree at most $2 \times (3n-6) + 1 = 6n - 11$, with finite ramification points at the
vertices of $P$ and no others. If the projection mapping is denoted by $\pi : \cl S^* \to \dd C$, then on $D^*$ we
have the equality $\pi \circ \tilde{F} = F$.}

\stit{Proof :} The required log-Riemann surface $\cl S$ may be constructed as follows:

\medskip

Let the circle of prime ends $\cl E$ be the union of $N = 3n - 6$ arcs $\alpha_1, \dots, \alpha_N$ say, ordered
cyclically so that each $\alpha_j$ intersects only $\alpha_{j-1}$ and $\alpha_{j+1}$, at its endpoints denoted by
$w^*_j$ and $w^*_{j+1}$ respectively (here and in what follows, indices are taken cyclically modulo $N$). For
$j=1,\dots,n$, let $F$ map $\alpha_j$ isometrically to a compact Euclidean segment $[w_j, w_{j+1}] \subset \dd C$
contained in a straight line $l_j \subset \dd C$. The line $l_j$ determines two half-planes; of these two, there
is one such that any point in $\alpha_j - \{ w^*_{j}, w^*_{j+1} \}$ has a neighbourhood in $D^*$ whose image under
$F$ is contained in the closure of this half-plane. We consider the {\it other} half-plane, and denote it by
$H_j$.

\medskip

For $j = 1, \dots, N$ we let $\theta_j, \phi_j \in [0, 2\pi)$ be angles such that $$\theta_j = \arg (w_{j+1} -
w_j) \ , \ \phi_j = \arg(w_{j-1} - w_j). $$ The half-plane $H_j$ is given by either $H_j = \{ \theta_j < \arg(w -
w_j) < \theta_j + \pi \}$ or $H_j = \{ \theta_j - \pi < \arg(w - w_j) < \theta_j \}$. If the first case occurs
then we let $m_j \in \dd Z$ be the smallest integer such that $2m_j \pi > \theta_j + \pi$ and define the angular
sector $U_j \subset \dd C$ by $U_j := \{ \theta_j + \pi < \arg(w - w_j) < 2m_j \pi \}$, while if the second case
occurs we take $m_j$ to be the largest integer such that $2 m_j \pi < \theta_j - \pi$ and define $U_j := \{
 2m_j \pi < \arg(w - w_j) < \theta_j - \pi \}$.

\medskip

Similarly, the half-plane $H_{j-1}$ is given by either $H_{j-1} = \{ \phi_j < \arg(w - w_j) < \phi_j + \pi \}$ or
$H_{j-1} = \{ \phi_j - \pi < \arg(w - w_j) < \phi_j \}$. In the first case we take $n_j \in \dd Z$ to be the
smallest integer such that $2 n_j \pi > \phi_j + \pi$ and define the angular sector $V_j \subset \dd C$ by $V_j :=
\{ \phi_j + \pi < \arg(w - w_j) < 2 n_j \pi \}$, while in the second case we take $n_j$ to be the largest integer
such that $2 n_j \pi < \phi_j - \pi$ and define $V_j := \{ 2n_j \pi < \arg(w - w_j) < \phi_j - \pi \}$.

\medskip

We paste $\overline{U_j} - \{ w_j \}$ to $\overline{H_j} - \{ w_j, w_{j+1} \}$ isometrically by the identity along
the open half-line common to their boundaries, which is either $\{ \arg(w - w_j) = \theta_j + \pi \}$ or $\{
\arg(w - w_j) = \theta_j - \pi \}$. Similarly we paste $\overline{V_j} - \{ w_j \}$ to $\overline{H_{j-1}}$
isometrically by the identity along either $\{ \arg(w - w_j) = \phi_j + \pi \}$ or $\{ \arg(w - w_j) = \phi_j -
\pi \}$. We paste $\overline{U_j} - \{ w_j \}$ to $\overline{V_j} - \{ w_j \}$ isometrically by the identity along
the remaining boundary segments $\{ \arg(w - w_j) = 2 m_j \pi \}$ and $\{ \arg(w - w_j) = 2 n_j \pi \}$. Finally
we paste $D^* - \{ w^*_1, \dots, w^*_n \}$ to each half-plane $\overline{H_j} - \{ w_j, w_{j+1} \}$, by
identifying the boundary segments $\alpha_j - \{ w^*_j, w^*_{j+1} \}$ and $(w_j, w_{j+1})$ isometrically via the
map $F$.

\medskip

The result is a space

$$\eqalign{
 \cl S & = (D^* - \{ w^*_1, \dots, w^*_N \}) \cr & \sqcup (\overline{H_1} - \{ w_1, w_2 \}) \sqcup \dots \sqcup
(\overline{H_N} - \{ w_N, w_1 \}) \cr & \sqcup (\overline{U_1} - \{ w_1 \}) \sqcup \dots \sqcup (\overline{U_N} -
\{ w_N \}) \cr & \sqcup (\overline{V_1} - \{ w_1 \}) \sqcup \dots \sqcup (\overline{V_N} - \{ w_N \}) / \sim \cr }
$$

where the equivalence relation $\sim$ describes how the various boundary segments are pasted.

\medskip

Since all the pasting maps are conformal, the space $\cl S$ inherits a Riemann surface structure from its
constituent parts. Moreover, each constituent comes with a distinguished chart, namely the identity map on the
domains $\overline{H_j} - \{ w_j, w_{j+1} \}, \overline{U_j} - \{ w_j \}, \overline{V_j} - \{ w_j \},
j=1,\dots,n$, and the map $F$ on $D^* - \{ w^*_1, \dots, w^*_N \}$. It is clear that these maps paste together to
give a globally defined map $\pi : \cl S \to \dd C$, which is a local diffeomorphism. If we denote by $\tilde{F} :
D^* - \{ w^*1, \dots, w^*_N \} \to \cl S$ the canonical embedding of $D^* - \{ w^*_1, \dots , w^*_N \}$ in $\cl
S$, then on $D$ we have $\pi \circ \tilde{F} = F$. We denote by $P$ the isometric image of $D$ in $\cl S$, $P =
\tilde{F}(D) \subset \cl S$.

\medskip

Pulling back the flat metric on $\dd C$ via $\pi$ gives a flat metric on $\cl S$, which in turn induces a global
metric on $\cl S$. The completion $\cl S^*$ of $\cl S$ with respect to this metric is given by adding $N$ points
$p_1, \dots, p_N$, $\cl S^* = \cl S \sqcup \{ p_1, \dots, p_N \}$, and each point $p_j$ is given by $p_j =
\partial H_j \cap \partial U_j \cap \partial V_j \cap \partial H_{j-1} \cap \partial P$ (with boundaries here
taken in $\cl S^*$. The embedding $\tilde{F} : D^* - \{ w^*1, \dots, w^*_N \} \to \cl S$ extends to an isometric
embedding $\tilde{F} : D^* \to P \cup \partial P \subset \cl S$ mapping $w^*_j$ to $p_j$, $j=1,\dots,N$, and the
equality $\pi \circ \tilde{F} = F$ continues to hold on $D^*$.

\medskip

Using the map $\pi$ and the metric on $\cl S$ it is shown in [Bi-PM] how to construct {\it log-charts } for $\cl
S$ and give $\cl S$ a compatible log-Riemann surface structure, for which $\cl S$ has $\pi$ as projection mapping.
The domain $P \subset \cl S$ is simply connected and its boundary in $\cl S^*$ is a finite union of compact
Euclidean segments, hence $P$ is a log-polygon; the vertices of $P$ are the points $p_1, \dots p_N$.

\medskip

Finally it is not hard to see that each point $p_j$ is a finite ramification point, and its order $k_j$ is at most
equal to $3$. As is explained in [Bi-PM], for any log-Riemann surface, the topological surface obtained by adding
all of its finite ramification points inherits a unique Riemann surface structure compatible with that of the
original log-Riemann surface. Thus in this case $\cl S^*$ inherits a Riemann surface structure from $\cl S$.
Moreover the surface $\cl S^*$ is simply connected; since $\cl S$ has a finite number of ramification points all
of finite order, it follows from a result in [Bi-PM] that $\cl S$ is the log-Riemann surface of a polynomial.
Moreover the degree of this polynomial is given by the sum $(k_1 - 1) + \dots + (k_N - 1) + 1$ which is at most
equal to $2N + 1$. $\diamondsuit$

\bigskip

{\bf 5.2.4) Construction of $\cl S_R$ and proof of the Main Theorem.}

\medskip

We can now construct the tube-log Riemann surface $\cl S_R$ and in the process prove the Main Theorem as follows:

\medskip

\stit{Proof of Main Theorem:} Let $\overline{C_1}, \dots, \overline{C_n}$ be closed half-cylinders isometric to
the pole-petals $\overline{P_1}, \dots, \overline{P_n}$ respectively, with $\overline{C_j} = \{ \hbox{ Re } (w /
\lambda_j) \leq 0 \} \subset \dd C / \dd 2\pi i \lambda_j \dd Z$. Let $\overline{P} = P \cup \partial P$ be the
closed log-polygon with projection $\pi : \overline{P} \to \dd C$ constructed in the previous section, and
$\tilde{F} : D \to P$ the isometry mapping $D$ to $P$ such that $F = \pi \circ \tilde{F}$, where $F$ is the
primitive defined in the previous section.

\medskip

Each petal boundary $\partial P_j$ corresponds to a side of $P$ isometric to a Euclidean segment with endpoints
differing by $2\pi i \lambda_j$, which we denote by $\beta_j$ say, $j = 1, \dots, n$. Each curve $\gamma_k$
corresponds to two sides of $P$ both isometric to $\gamma_k$, which we denote by $\Gamma^+_k, \Gamma^-_k$ say, $k
= 1, \dots, n-3$.

\medskip

For $j=1, \dots, n$, let $\overline{C_j}' := \overline{C_j} - \{ a_j \}$ be obtained from $\overline{C_j} \subset
\dd C / \dd 2\pi i \lambda_j \dd Z$ by deleting a point $a_j \in \partial C_j \subset \dd C / \dd 2\pi i \lambda_j
\dd Z$.  Let $\overline{P}' := \overline{P} - \{ p_1, \dots, p_N \}$ be obtained from $\overline{P}$ by deleting
the vertices of $\overline{P}$.

\medskip

For $j = 1,\dots,n$ the open Euclidean segments $\partial C_j - \{ a_j \}$, $\beta_j \cap \overline{P}'$ of
$\overline{C_j}'$ and $\overline{P}'$ respectively are isometric, and we paste them together isometrically along
these segments. For $k = 1, \dots, n-3$ we paste the open Euclidean segments $\Gamma^+_k \cap \overline{P}'$ and
$\Gamma^-_k \cap \overline{P}'$ of $\overline{P}'$ isometrically to each other. The resulting space $$ \cl S :=
\overline{P}' \sqcup (\overline{C_1}' \sqcup \dots \sqcup \overline{C_n}') / \sim $$ (where the relation $\sim$
denotes the various isometric pastings) is clearly a tube-log Riemann surface. Moreover it is the tube-log Riemann
surface associated to $\int R(z) \ dz$, $\cl S = \cl S_R$; the required biholomorphism $\tilde{F} : \overline{\dd
C} - (Z \cup P) \to \cl S$ may be constructed as follows:

\medskip

On $D$ we let $\tilde{F}$ be given by the isometry (already denoted by $\tilde{F}$) mapping $D$ to the log-polygon
$P$, $\tilde{F} : D \to P$. For $j=1, \dots , n$, this map extends continuously to every point of $\partial P_j$
which is not a critical point of $R(z) \ dz$, mapping $\partial P_j - Z$ isometrically to $\beta_j \cap
\overline{P}'$ For $k=1, \dots, n-3$, at each point $z_0 \in \gamma_k$ which is not an endpoint of $\gamma_k$,
there are two prime ends of $D$ with impression $z_0$; when $z \to z_0$ through one of these prime ends,
$\tilde{F}(z)$ tends to a point of $\Gamma^+_k \cap \overline{P}'$, and when $z \to z_0$ through the other, to a
point of $\Gamma^-_k \cap \overline{P}'$. These two points correspond under the isometric pasting of $\Gamma^+_k
\cap \overline{P}'$ and $\Gamma^-_k \cap \overline{P}'$. It follows from these remarks that $\tilde{F} : D \to P$
extends to an isometric embedding into $\cl S$,  $\tilde{F} : D \cup ((\partial P_1 \cup \dots \partial P_n) - Z)
\cup(\gamma_1 \cup \dots \cup \gamma_{n-3}) \to (\overline{P}' / \sim) \subset \cl S$ (isometric for $D \cup
(\gamma_1 \cup \dots \cup \gamma_{n-3})$ considered as a subspace of $(X, d)$).

\medskip

For $j = 1, \dots, n$, the closed pole-petal $\overline{P_j}$ (with the metric induced by $|R(z)| |dz|$) and the
half-cylinder $\overline{C_j}'$ are isometric, and moreover, removing the critical point which lies on $\partial
P_j$, we can choose an isometry $\tilde{F_j} : \overline{P_j} - Z \to \overline{C_j}'$ which agrees with the map
$\tilde{F}$ above on $\partial P_j - Z$.

\medskip

Putting together the maps $\tilde{F}$ and $\tilde{F_j}, j=1, \dots,n$ gives the required biholomorphism, which we
denote by $\tilde{F} : \overline{\dd C} - (Z \cup P) \to \cl S$. The derivative $\tilde{F}' : \overline{\dd C} -
(Z \cup P) \to \dd C$ computed in the distinguished charts on $\cl S$ is clearly equal to the rational function
$R(z)$ on $D$, and hence on all of $\overline{\dd C} - (Z \cup P)$. Thus $\cl S = \cl S_R$. Moreover it is clear
that the tube-log Riemann surface $\cl S_R$ is constructed as stated in the Main Theorem. $\diamondsuit$

\bigskip

\centerline{{\bf References}}

\bigskip

[Bi-PM] K.BISWAS, R.PEREZ-MARCO, {\it Log-Riemann Surfaces}, Preprint, UCLA, 2004.

\medskip

[PM1] R.PEREZ-MARCO,  {\it Sur les dynamiques holomorphes non lin\'earisables et une conjecture de V.I. Arnold},
Ann. Scient. Ec. Norm. Sup. 4 serie, {\bf 26}, 1993, p.565-644; (C.R. Acad. Sci. Paris, {\bf 312}, 1991,
p.105-121).

\medskip

[PM2] R.PEREZ-MARCO, {\it Uncountable number of symmetries for non-linearisable holomorphic dynamics}, Inventiones
Mathematicae, {\bf 119}, {\bf 1}, p.67-127,1995; (C.R. Acad. Sci. Paris, {\bf 313}, 1991, p. 461-464).

\medskip

[PM3] R.PEREZ-MARCO, {\it Siegel disks with smooth boundary}, Preprint, Universit\'e de PARIS-SUD, 1997,
www.math.ucla.edu/$\tilde{\ }$ricardo, submitted to the Annals of Mathematics in 1997.

\end